%




\input amstex.tex

\magnification=\magstep1
\hsize=5.5truein
\vsize=9truein
\hoffset=0.5truein
\parindent=10pt
\newdimen\nagykoz
\newdimen\kiskoz
\nagykoz=7pt
\kiskoz=2pt
\parskip=\nagykoz
\baselineskip=12.7pt


\loadeufm \loadmsam \loadmsbm

\font\vastag=cmssbx10
\font\drot=cmssdc10
\font\vekony=cmss10
\font\vekonydolt=cmssi10
\font\cimbetu=cmssbx10 scaled \magstep1
\font\szerzobetu=cmss10

\font\scVIII=cmcsc8
\font\rmVIII=cmr8
\font\itVIII=cmti8
\font\bfVIII=cmbx8
\font\ttVIII=cmtt8

\def\cim#1{{\centerline{\cimbetu#1}}}
\def\szerzo#1{{\vskip0.3truein\centerline{\szerzobetu#1}}}
\def\alcim#1{{\medskip\centerline{\vastag#1}}}
\def\tetel#1#2{{{\drot#1}{\it\szukebb~#2\tagabb}}}
\long\def\biz#1#2{{{\vekony#1} #2}}
\def\kiemel#1{{\vekonydolt#1\/}}
\long\def\absztrakt#1#2{{\vskip0.4truein{\vekony#1} #2\vskip0.5truein}}
\def\szukebb{\parskip=\kiskoz}
\def\tagabb{\parskip=\nagykoz}
\def\vonal{{\vrule height 0.2pt depth 0.2pt width 0.5truein}}

\def\CC{{\Bbb C}}
\def\Cb{{C_{\text b}}}
\def\Cu{{C_{\text u}}}
\def\Cbu{{C_{\text{bu}}}}
\def\OObu{{\OO_{\text{bu}}}}
\def\Cs{C_{\text s}}
\def\Zs{Z_{\text s}}
\def\Hs{H_{\text s}}

\def\ts{\textstyle}

\def\bmfd{{Banach manifold}}
\def\balg{{Banach algebra}}
\def\blgp{{Banach Lie group}}
\def\bla{{Banach Lie algebra}}
\def\bvbdl{{Banach vector bundle}}
\def\hvbdl{{Hilbert vector bundle}}
\def\blgpbdl{{Banach Lie group bundle}}
\def\blabdl{{Banach Lie algebra bundle}}

\def\re{\hbox{\rm Re}}
\def\im{\hbox{\rm Im}}

\def\bigO{{\text{\rm O}}}

\def\Oa{{\Omega}}
\def\oa{{\omega}}
\def\UU{{\frak U}}
\def\VV{{\frak V}}
\def\BB{{\frak B}}
\def\GL{\hbox{\rm GL}}
\def\diag{\text{\rm diag}}

\def\Hom{\hbox{\rm Hom}}
\def\End{\hbox{\rm End}}

\def\Ad{\hbox{\rm Ad}}

\def\vbdl{{vector bundle}}

\def\cts{{continuous}}
\def\ucts{{uniformly continuous}}
\def\bu{{bounded and uniformly continuous}}
\def\idml{{infinite dimensional}}

\def\idms{{infinite dimensions}}

\def\fdml{{finite dimensional}}

\def\fdms{{finite dimensions}}

\def\pscx{{pseudoconvex}}

\def\dist{{\hbox{\rm dist}}}

\def\sbs{{Schauder basis}}

\def\cubs{{countable unconditional basis}}
\def\psh{{plurisubharmonic}}
\def\pshdom{{\psh\ domination}}

\def\st{{such that}}
\def\wrt{{with respect to}}
\def\subsub{\subset\!\subset}
\def\delbar{{\bold{\bar\partial}}}
\def\llip{{C^{0,1}_{\text{loc}}}}

\def\<{{\langle}}
\def\>{{\rangle}}
\def\NN {{\Bbb N}}
\def\RR {{\Bbb R}}

\def\AA{{\Cal A}}

\def\OO {{\Cal O}}
\def\OOs{{\Cal O}_{\text s}}

\def\ho{{(h${}_1$)}}

\def\ro{\varrho}

\def\epsz{\varepsilon}
\def\fii{\varphi}
\def\fn{func\-tion}
\def\fns{func\-tions}
\def\holo{hol\-o\-mor\-phic}
\def\tpl{topological}

\def\cpx{complex}
\def\ses{short exact sequence}

\def\nbd{neighbor\-hood}

\def\bspc{Banach space}

\def\ka{{\kappa}}
\def\la{\lambda}
\def\La{\Lambda}
\def\sa{\sigma}

\def\da{{\delta}}

\def\ga{{\gamma}}
\def\aa{{\alpha}}
\def\ba{{\beta}}
\def\da{{\delta}}
\def\ta{{\theta}}

\def\Ga{{\Gamma}}

\def\za{\zeta}

\def\Prop{Proposition}
\def\p#1{{\Prop~#1}}

\def\Th{{Theorem}}

\def\t#1{{\Th~#1}}

{\phantom.}
\vskip0.5truein
\cim{ON HOLOMORPHIC BANACH VECTOR}
\vskip0.2truein
\cim{BUNDLES OVER BANACH SPACES}
\szerzo{Imre Patyi\plainfootnote{${}^*$}{\rmVIII Supported in part by a Research Initiation Grant
	from Georgia State University.}}
\absztrakt{ABSTRACT.}{Let $X$ be a \bspc\ with a \cubs\ (e.g., $X=\ell_2$ Hilbert space),
	$\Oa\subset X$
	\pscx\ open, $E\to\Oa$ a locally trivial \holo\ \vbdl\ with a
	\bspc\ $Z$ for fiber type, $\OO^E$ the sheaf of germs of
	\holo\ sections of $E\to\Oa$, and
	$Z_1$ the \bspc\
	$Z_1=\ell_p(Z)=
	\{z=(z_n):z_n\in Z,\ \|z\|=\big(\sum_{n=1}^\infty\|z_n\|^p\big)^{1/p}<\infty\}$, $1\le p<\infty$.
	 Then $E\oplus(\Oa\times Z_1)$ and $\Oa\times Z_1$ are \holo{}ally isomorphic,
	$\OO^E$ is acyclic and $E$ is so to speak stably trivial over $\Oa$ in a generalized sense.
	 We also show that if $E$ is \cts{}ly trivial over $\Oa$, then
	$E$ is \holo{}ally trivial over $\Oa$.
	 In particular, if $Z=\ell_2$ or $\Oa$ is contractible, then $E$ is 
	\holo{}ally trivial over $\Oa$.
	 Some applications are also given.

	 MSC 2000: 32L05, (32L10, 32L20, 32Q28, 46G20)
	
	 Key words: analytic cohomology, pseudoconvex domains, holomorphic Banach vector bundles.


\hfill{\it To my dear little daughter S\'ari Mangala on her third birthday.}
}

\def\sA{{1}}
\def\sB{{2}}
\def\sC{{3}}
\def\sD{{4}}
\def\sE{{5}}
\def\sF{{6}}
\def\sG{{7}}
\def\sH{{8}}
\def\sI{{9}}
\def\sJ{{10}}
\def\sK{{11}}
\def\sL{{12}}
\def\sM{{14}}
\def\sN{{15}}
\def\sX{{13}}


\def\tAA{{\sA.1}}
\def\tAB{{\sA.2}}
\def\tAC{{\sA.3}}

\def\tBA{{\sB.1}}
\def\tBB{{\sB.2}}

\def\eBA{{\sB.1}}
\def\eBB{{\sB.2}}

\def\tCA{{\sC.1}}
\def\tCB{{\sC.2}}
\def\tCC{{\sC.3}}

\def\eCA{{\sC.1}}

\def\tDA{{\sD.1}}

\def\tEA{{\sE.1}}
\def\tEB{{\sE.2}}
\def\tEC{{\sE.3}}
\def\tED{{\sE.4}}
\def\tEE{{\sE.5}}
\def\tEF{{\sE.6}}
\def\tEG{{\sE.7}}
\def\tEH{{\sE.8}}
\def\tEI{{\sE.9}}

\def\eEA{{\sE.1}}
\def\eEB{{\sE.2}}
\def\eEC{{\sE.3}}
\def\eED{{\sE.4}}
\def\eEE{{\sE.5}}

\def\tFA{{\sF.1}}
\def\tFB{{\sF.2}}
\def\tFC{{\sF.3}}
\def\tFD{{\sF.4}}
\def\tFE{{\sF.5}}
\def\tFF{{\sF.6}}
\def\tFG{{\sF.7}}
\def\tFH{{\sF.8}}

\def\tGA{{\sG.1}}
\def\tGB{{\sG.2}}
\def\tGC{{\sG.3}}
\def\tGD{{\sG.4}}

\def\tHA{{\sH.1}}
\def\tHB{{\sH.2}}
\def\tHC{{\sH.3}}
\def\tHD{{\sH.4}}
\def\tHE{{\sH.5}}
\def\tHF{{\sH.6}}
\def\tHG{{\sH.7}}
\def\tHH{{\sH.8}}
\def\tHI{{\sH.9}}
\def\tHJ{{\sH.10}}

\def\tHA{{\sH.1}}
\def\tHB{{\sH.2}}
\def\tHC{{\sH.3}}
\def\tHD{{\sH.4}}
\def\tHE{{\sH.5}}

\def\eHA{\sH.1}
\def\eHB{\sH.2}
\def\eHC{\sH.3}

\def\tLA{{\sL.1}}

\def\tMA{{\sM.1}}
\def\tMB{{\sM.3}}
\def\tMC{{\sM.4}}
\def\tMD{{\sM.5}}
\def\tME{{\sM.6}}
\def\tMF{{\sM.7}}
\def\tMG{{\sM.8}}
\def\tMH{{\sM.9}}
\def\tMI{{\sM.10}}
\def\tMX{{\sM.2}}

\def\rC{C}
\def\rH{H}
\def\rK{K}
\def\rLt{Lt}
\def\rLA{L1}
\def\rLB{L2}
\def\rLC{L3}
\def\rLD{L4}
\def\rN{N}
\def\rPA{P1}
\def\rPB{P4}
\def\rPC{P5}
\def\rPD{P6}
\def\rPE{P3}
\def\rPF{P2}
\def\rZ{F}

\alcim{\sA. INTRODUCTION.}

	 In this paper we study \holo\ \bvbdl{}s over \pscx\ open subsets of \bspc{}s.
	 Under suitable conditions we show that any \holo\ \bvbdl\ can be exhibited as
	 a direct summand
	of a trivial \holo\ \bvbdl.
	 In this way we recover and substantially extend the following theorem of Lempert.

\tetel{\t\tAA.}{{\rm(Lempert, [\rLC])}
	Let $X$ be a \bspc\ with a \cubs, $\Oa\subset X$ \pscx\ open, $E\to\Oa$
	a \holo\ \bvbdl, $\OO^E\to\Oa$ the sheaf of germs of \holo\ sections of $E\to\Oa$.
	Then the sheaf cohomology groups $H^q(\Oa,\OO^E)$ vanish for all $q\ge1$.
}

	 Following [\rLB] by Lempert we say that \kiemel{\pshdom} holds in a
	complex \bmfd\ $\Oa$ if for every $u:\Oa\to\RR$ locally upper bounded
	there is a $\psi:\Oa\to\RR$ \cts\ and \psh\ \st\ $u(x)<\psi(x)$ for all $x\in\Oa$.

\tetel{\t\tAB.}{{\rm(Lempert, [\rLB])}
	If $X,\Oa$ are as in \t\tAA, then \pshdom\ holds in\/ $\Oa$.
}

	 We prove here the following \t\tAC.

\tetel{\t\tAC.}{Let $X$ be a \bspc\ with a \sbs, $\Oa\subset X$ \pscx\ open,
	$E\to\Oa$ a \holo\ \bvbdl\ with a \bspc\ $Z$ for fiber type.
	If \pshdom\ holds in\/ $\Oa$, then we have the following.
\vskip0pt
	{\rm(a)} $H^q(\Oa,\OO^Z)=0$ for $q\ge1$.
\vskip0pt
	{\rm(b)} Let $Z_1=\ell_p(Z)$, $1\le p<\infty$.
	Then $E\oplus(\Oa\times Z_1)$ and\/ $\Oa\times Z_1$ are
	\holo{}ally isomorphic over\/ $\Oa$.
\vskip0pt
	{\rm(c)} $H^q(\Oa,\OO^E)=0$ for $q\ge1$.
\vskip0pt
	{\rm(d)} If $E$ is \cts{}ly trivial over\/ $\Oa$, then
	$E$ is \holo{}ally trivial over\/ $\Oa$.
\vskip0pt
	{\rm(e)} If $Z=\ell_2$, then $E$ is \holo{}ally trivial over\/ $\Oa$.
\vskip0pt
	{\rm(f)} If\/ $\Oa$ is contractible, then $E$ is \holo{}ally
	trivial over\/ $\Oa$.
}

	 The proof of \t\tAC\ depends on an exhaustion procedure by Lempert,
	on a simple aspect of the Grauert--Oka
	principle in [\rPA], and on the contractibility of $\GL(\ell_2)$
	by Kuiper in [\rK].
	 A few applications of \t\tAC\ are given at the end of this paper.
	 The methods and results here are further applied in [\rPB] on an
	analytic Koszul complex, and in [\rPC] on sheaf
	cohomology vanishing for a general class of sheaves that includes
	vector bundles and certain ideal sheaves.
	 For background see [\rLA--\rLD, \rPA].

\alcim{\sB. EXHAUSTION.}

	 This section describes a way to exhaust a \pscx\ open subset $\Oa$ of a \bspc\ $X$
	that is convenient for proving vanishing results for sheaf cohomology over $\Oa$.
	 We follow here [\rLD, \S\,2].

	 We say that a \fn\ $\aa$, call their set $\AA'$, is an
	\kiemel{admissible radius \fn} on $\Oa$ if $\aa:\Oa\to(0,1)$ is \cts\
	and $\aa(x)<\dist(x,X\setminus\Oa)$ for $x\in\Oa$.
	 We say that a \fn\ $\aa$, call their set $\AA$, is an
	\kiemel{admissible Hartogs radius \fn} on $\Oa$ if $\aa\in\AA'$
	and $-\log\aa$ is \psh\ on $\Oa$.
	 Call $\AA$ \kiemel{cofinal} in $\AA'$ if for each $\aa\in\AA'$
	there is a $\ba\in\AA$ with $\ba(x)<\aa(x)$ for $x\in\Oa$.

\tetel{\p\tBA.}{Plurisubharmonic domination holds in\/ $\Oa$ if and only if $\AA$
	is cofinal in $\AA'$.
}

\biz{Proof.}{Write $\aa=e^{-u}\in\AA'$ and $\ba=e^{-\psi}\in\AA$.
	As \pshdom\ holds on $\Oa$ for $u$ \cts\ if and only if for
	$u$ locally upper bounded, the proof of \p\tCA\ is complete.
}

	 It will be often useful to look at coverings by balls
	$B_X(x,\aa(x))$, $x\in\Oa$, $\aa\in\AA'$, and shrink their
	radii to obtain a finer covering by balls $B_X(x,\ba(x))$,
	$x\in\Oa$, $\ba\in\AA$.

	 Let $e_n$, $n\ge1$, be a \sbs\ in the \bspc\ $(X,\|\cdot\|)$.
	 One can change the norm $\|\cdot\|$ to an equivalent
	norm so that
$
	\big\|\sum_{i=m}^nx_ie_i\big\|\le\big\|\sum_{i=M}^Nx_ie_i\big\|
 $
 	for $0\le M\le m\le n\le N\le\infty$, $x_i\in\CC$.
	 Introduce the projections $\pi_N:X\to X$,
	$\pi_N\sum_{i=1}^\infty x_ie_i=\sum_{i=1}^Nx_ie_i$, 
	$x_i\in\CC$, $\pi_0=0$, $\pi_\infty=1$, $\ro_N=1-\pi_N$, and define for
	$\aa\in\AA$ and $N\ge0$ integer the sets
$$\eqalign{
	D_N\<\aa\>&=\{\xi\in\Oa\cap\pi_N X:(N+1)\aa(\xi)>1\},\cr
	\Oa_N\<\aa\>&=\{x\in\pi_N^{-1}D_N\<\aa\>:\|\ro_N x\|<\aa(\pi_N x)\},\cr
	D^N\<\aa\>&=\pi_{N+1}X\cap\Oa_N\<\aa\>,\cr
	\Oa^N\<\aa\>&=\{x\in\pi_{N+1}^{-1}D^N\<\aa\>:\|\ro_{N+1}x\|<\aa(\pi_Nx)\},\cr
	\BB(\aa)&=\{B_X(x,\aa(x)):x\in\Oa\},\cr
	\BB_N(\aa)&=\{B_X(x,\aa(x)):x\in\Oa_N\<\aa\>\}.\cr
}
\tag\eBA
 $$

	 These $\Oa_N\<\aa\>$ are \pscx\ open in $\Oa$, and they will serve
	to exhaust $\Oa$ as $N=0,1,2,\ldots$ varies.

\tetel{\p\tBB.}{{\rm(Lempert)} Let $\aa\in\AA$, and suppose that \pshdom\ holds in\/ $\Oa$.
\vskip0pt  
	{\rm(a)} There is an $\aa'\in\AA$, $\aa'<\aa$, with\/
	$\Oa_n\<\aa'\>\subset\Oa_N\<\aa\>$ for all $N\ge n$.
	 So any $x_0\in\Oa$ has a \nbd\ contained in all but
	finitely many $\Oa_N\<\aa\>$.
\vskip0pt  
	{\rm(b)} There are $\ba,\ga\in\AA$, $\ga<\ba<\aa$, so that
	for all $N$ and $x\in\Oa_N\<\ga\>$
$$
	B_X(x,\ga(x))\subset\Oa_N\<\ba\>\cap\pi_N^{-1}B_X(\pi_N x,\ba(x))
	\subset B_X(x,\aa(x)).
\tag\eBB
 $$
\vskip0pt  
	{\rm(c)} If\/ $8\aa\in\AA$, $Y\subset X$ is a \fdml\ complex affine
	subspace, then $Y\cap\overline{\Oa_N\<\aa\>}$ is \psh{}ally convex
	in $Y\cap\Oa$.
\vskip0pt
	{\rm(d)} We have that\/ $\overline{\Oa_N\<\aa\>}\subset\overline{\Oa^N\<\aa\>}$.
	If\/ $4\aa\in\AA$, then\/ $\Oa^N\<\aa\>\subset\Oa_N\<2\aa\>$.
\vskip0pt
	{\rm(e)} There is a $\ba\in\AA$, $\ba<\aa$, with
	$\Oa_N\<\ba\>\subset\Oa_N\<\aa\>\cap\Oa_{N+1}\<\aa\>$ for $N\ge0$.
\vskip0pt
	{\rm(f)} There is an $\aa'\in\AA$, $\aa'<\aa$, \st\
	the covering $\BB_N(\aa)|\Oa_N\<\aa'\>$ has a finite basic refinement
	for all $N\ge0$.
}

\biz{Proof.}{For (a) and (b) see [\rLD, Prop.\,2.1], and [\rLC, Prop.\,4.3],
	for (c) [\rLC, Prop.\,4.3], for (d) [\rLC, Prop.\,4.4], for (e) [\rLD, Prop.\,2.3],
	and for (f) see [\rPB, Prop.\,3.2(c)].  
	 The proof of \p\tBB\ is complete.
	 (Remark for the record that (f) was not explicitly formulated by Lempert.)
}

	The meaning of \p\tBB(b) is that certain refinement maps exist
	between certain open coverings, while (cd) are useful for Runge
	type approximation, and (ef) for exhaustion.

\alcim{\sC. APPROXIMATION.}

	This section deals with some versions of Runge approximation.

\tetel{\p\tCA.}{Let $X,Z$ be \bspc{}s, $f\in\OO(B_X(1),Z)$ bounded, $0<\ta<1$, $\epsz>0$.
	Then there is a polynomial $g\in\OO(X,Z)$ with\/ $\|f(x)-g(x)\|<\epsz$
	for\/ $\|x\|<\ta$.
}

\biz{Proof.}{As well known we may take $g=\sum_{m=0}^n f_m$ for $n$ large enough,
	where $\sum_{m=0}^\infty f_m$ is the homogeneous expansion of $f$ about $x=0$.
}

	Let $Y,Z$ be \bspc{}s, $D_1\subsub D_2\subsub D_3\subsub\CC^N$ \pscx\ open,
	$R:\overline{D_3}\to(0,\infty)$ \cts\ with $-\log R$ \psh\ on $D_3$.
	Let
$$
	\Oa(D_3,R)=\{(\za,y)\in D_3\times Y:\|y\|<R(\za)\}.
\tag\eCA
 $$

\tetel{\p\tCB.}{If $\overline{D_2}$ is \holo{}ally convex in $D_3$, $f\in\OO(\Oa(D_2,R),Z)$ is
	\bu, $0<\ta<1$, and $\epsz>0$, then there is a $g\in\OO(D_3\times Y,Z)$ that is \bu\ on
	any set of the form $K\times B_Y(r)$, where\/ $0<r<\infty$ and $K\subsub D_3$,
	and\/ $\|f(\za,y)-g(\za,y)\|<\epsz$ for $(\za,y)\in\Oa(D_1,\ta R)$.
}

\biz{Proof.}{See the proof of [\rLA, Thm.\,6.1] and replace the reference there
	to Hypothesis$(X,F)$ by reference to \p\tCA\ here.
}

	Let $X$ be as in \S\,\sB.

\tetel{\p\tCC.}{Let\/ $8\aa\in\AA$, and choose $\ga\in\AA$ as in \p\tBB(b), and $\ga'\in\AA$
	as in \p\tBB(a).
	 Let $f\in\OO(\Oa_N\<\aa\>,Z)$ be \bu, and $\epsz>0$.
\vskip0pt
	{\rm(a)} There is a $g\in\OO(\Oa_{N+1}\<\aa\>,Z)$ \bu\ with\/
	$\|f(x)-g(x)\|<\epsz$ for $x\in\Oa_N\<\ga\>$.
\vskip0pt
	{\rm(b)} There is a $g\in\OO(\Oa,Z)$ that is \bu\ on\/ $\Oa_{N+p}\<\ga'\>$
	for all $p\ge0$ with\/ $\|f(x)-g(x)\|<\epsz$ for $x\in\Oa_N\<\ga'\>$.
}

\biz{Proof.}{(a) Apply \p\tBB\ and \p\tCB\ as in the proof of [\rLC, Thm.\,4.5].

	(b) Letting $g_0=f$, apply (a) repeatedly to find \bu\ \fns\ $g_p\in\OO(\Oa_{N+p}\<\aa\>,Z)$
	for $p\ge1$ with $\|g_{p}(x)-g_{p-1}(x)\|<\epsz/2^{p+1}$ for $x\in\Oa_{N+p}\<\ga\>$.
	Letting $g=g_0+\sum_{p=1}^\infty(g_{p}-g_{p-1})=\lim_{p\to\infty}g_p$
	completes the proof of \p\tCC.
}

\alcim{\sD. VANISHING IN THE MIDRANGE.}

	In this section we show that certain cocycles can be resolved.

	Let $D\subsub D'\subsub\CC^N$ be \pscx\ open, $\UU$ a covering of $D$
	by open sets $U\subset D$.
	 For $U\in\UU$ let $U'$ be \pscx\ open with $\overline U\subset U'\subset D'$
	and assume that $\UU'=\{U':U\in\UU\}$ covers $D'$.
	 Let $Y$ be a \bspc, $\pi:\CC^N\times Y\to\CC^N\times Y$ the projection
	$\pi(\za,y)=(\za,0)$, $R:D'\to(0,\infty)$ be \cts\ and bounded away from zero
	and $-\log R$ \psh, $0<\ta<1$, $\Oa'=\Oa(D',R)$, $\Oa=\Oa(D,\ta R)$ as in $(\eCA)$,
	$\UU'(\Oa')=\{U'(\Oa')=\pi^{-1}(U')\cap\Oa':U'\in\UU'\}$,
	$\UU(\Oa)=\{U(\Oa)=\pi^{-1}(U)\cap\Oa:U\in\UU\}$.

	 We say that a cochain $f\in C^q(\VV,\OO^Z)$, $q\ge0$, is \bu\ if
	each of its components $f_\sa$ is \bu\ on the body $|\sa|$ of any $q$-simplex
	$\sa$ of $\VV$.
	 This is most useful when the covering $\VV$ is finite or has a finite refinement.

\tetel{\p\tDA.}{For any \bu\ cocycle $f\in Z^q(\UU'(\Oa'),\OO^Z)$, $q\ge1$,
	there is a \bu\ cochain $g\in C^{q-1}(\UU(\Oa),\OO^Z)$ with $f|\UU(\Oa)=\da g$.
}

\biz{Proof.}{This can be done using a smooth partition of unity on $D'$ and
	an integral operator solving a $\delbar$-equation on a smooth strictly
	\pscx\ complete Hartogs domain.  See [\rPA, Prop.\,5.1, Prop.\,7.1].
}

\alcim{\sE. INVERTIBLE MATRICES.}

	This section describes some properties of \holo\ \fns\ with values in invertible
	matrices and invertible linear operators on a \bspc.

	Let $Z$ be a \bspc, $\End(Z)=\Hom(Z,Z)$ the \balg\ of all linear operators
	$A:Z\to Z$ endowed with the operator norm $\|A\|$, $G=\GL(Z)$ the \blgp\ of
	units of $\End(Z)$, and $\dot G=\End(Z)$ the \bla\ of $G$.
	 Here $G$ is a \pscx\ open subset of $\End(Z)$, and we have a \holo\ map,
	called the exponential map, $\exp:\dot G\to G$ defined by 
	$\exp(\xi)=\sum_{n=0}^\infty \xi^n/n!$, which is bi\holo\ from a small
	ball $B_{\dot G}(\epsz_0)$ of $\dot G$ to an open \nbd\ of $1\in G$, whose
	inverse is called logarithm and is denoted by $\log$.
	 In this case $\exp$ may equally be defined by 
	$\exp(\xi)=\lim_{n\to\infty}\big(1+\frac\xi n\big)^n$.
	
	 Let $(\Oa,d)$ be a metric space, and $f\in C(\Oa,G)$ a \cts\ \fn.
	If both $\|f(x)\|$ and $\|f(x)^{-1}\|$ are bounded for $x\in\Oa$,
	then we say that $f$ is
	\kiemel{bounded on} $\Oa$, and write that
	$f\in\Cb(\Oa,G)$.
	 If for any $\epsz>0$ there is a $\da>0$ \st\
	$\|f(x)f(y)^{-1}-1\|<\epsz$ for all $x,y\in\Oa$ with $d(x,y)<\da$,
	then we say that $f$ is
	\kiemel{\ucts\ on} $\Oa$, and write that
	$f\in\Cu(\Oa,G)$.
	 If $f$ is both \bu\ on $\Oa$, then we write that $f\in\Cbu(\Oa,G)$.
	 If $\Oa$ is a complex \bmfd\ with a fixed metric $d$ that induces the
	topology of $\Oa$, then we endow $\Oa\times[0,1]$ with the metric
	$\ro(x,t;y,s)=\max\{d(x,y),|t-s|\}$.
	 If $f\in\OO(\Oa,G)$ and there is an $h\in\Cbu(\Oa\times[0,1],G)$ \st\
	$h(x,1)=f(x)$ and $h(x,0)=1$ for all $x\in\Oa$, then we say that $f$ is a
	\kiemel{simple \fn\ on} $\Oa$ and that $h$ is a 
	\kiemel{simple \ho-\fn} (or a simple
	one-parameter homotopy)
	\kiemel{on} $\Oa$, and write that $f\in\OOs(\Oa,G)$ and that $h\in\OOs^1(\Oa,G)$.
	A simple \fn\ $f\in\OOs(\Oa,Z)$ and a simple \ho-\fn\ $h\in\OOs^1(\Oa,Z)$
	with values in a \bspc\ $Z$
	are similarly defined: $f$ is just a \holo\ \fn\ $f\in\OO(\Oa,Z)$ that is \bu\ on
	$\Oa$ (null homotopy is automatic by taking $tf(x)$), and $h$ is just a homotopy
	$h\in\Cbu(\Oa\times[0,1],Z)$ that is \holo\ in $x\in\Oa$ for $t=1$,
	and $h(x,0)=0$ for $x\in\Oa$.
	 Simple \fns\ form a group $\OOs(\Oa,G)$ \wrt\ pointwise multiplication and inversion
	that is closed under uniform sequential limits, and similarly for simple \ho-\fns\
	$\OOs^1(\Oa,G)$.
	 If $H$ is a topological group, then let $C_0([0,1],H)$ be the set of all \cts\ \fns\
	$f:[0,1]\to H$ with $f(0)=1\in H$.
 	 All that we state below in this section for $G=\GL(Z)$ can easily and quite
	analogously be extended to the \blgp\ $G=C_0([0,1],\GL(Z))$, the path group of
	$\GL(Z)$.

\tetel{\p\tEA.}{{\rm(a)} If $f,g\in\OOs(\Oa,G)$ are simple \fns, then so are $fg\in\OOs(\Oa,G)$
	and $f^{-1}\in\OOs(\Oa,G)$.
\vskip0pt
	{\rm(b)} If $f_n\in\OOs(\Oa,G)$ are simple \fns\ on\/ $\Oa$, and $f_n\to f$
	uniformly on\/ $\Oa$, i.e., $\lim_{n\to\infty}\sup_{x\in\Oa}\|f(x)f_n(x)^{-1}-1\|=0$,
	then their limit $f\in\OOs(\Oa,G)$ is also a simple \fn\ on\/ $\Oa$.
\vskip0pt
	{\rm(c)} If $f\in\OO(\Oa,G)$ is \st\ there is a homotopy $h\in\Cbu(\Oa\times[0,1],G)$
	with $h(x,t)$ \holo\ in $x\in\Oa$ for each $t\in[0,1]$, and $h(x,0)=1$ and $h(x,1)=f(x)$
	for $x\in\Oa$, then there are finitely many \fns\ $g_1,\ldots,g_n\in\OOs(\Oa,\dot G)$ \st\
	$f(x)=\exp(g_1(x))\ldots\exp(g_n(x))$ for all $x\in\Oa$.
\vskip0pt
	{\rm(d)} If $g\in\OOs^1(\Oa,\dot G)$, then the solution $h\in\OOs^1(\Oa,G)$
	of the parametric initial value problem $\frac{d}{dt}h(x,t)=g(x,t)h(x,t)$,
	$h(x,0)=1$ for the linear ordinary differential equation for the left
	logarithic derivative is a simple \ho-\fn. Thus $f(x)=h(x,1)$ defines a simple \fn\
	$f\in\OOs(\Oa,G)$.
}

\biz{Proof.}{As this proposition expresses some well known simple facts only, we just
	indicate the arguments.

	(a) Boundedness and null homotopy are clear.  Uniform continuity for $fg$
	follows on writing 
	$f(x)g(x)(f(y)g(y))^{-1}=[f(x)f(y)^{-1}][g(x)g(y)^{-1}]^{f(y)}$
	from considering conjugation and the adjoint action of $G$ on $\dot G$,
	where $a^b=bab^{-1}$.

	(b) The limit is clearly bounded and \ucts.
	To see that it is also null homotopic note that $\|f(x)f_n(x)^{-1}-1\|<1/2$
	for $x\in\Oa$ if $n$ is large enough.
	Fixing one such $n$ we can write $f(x)f_n(x)^{-1}=\exp(g(x))$,
	or $f(x)=\exp(g(x))f_n(x)$, from which $f$ is easily seen null homotopic.

	(c) To write $f$ as a finite product of exponentials of simple \fns\ we
	consider the telescopic product	
$$\eqalign{
	f(x)=h(x,{}&0)^{-1}h(x,1)=\cr
	&\ts[h(x,0)^{-1}h(x,\frac1n)][h(x,\frac1n)^{-1}h(x,\frac2n)]\ldots
	[h(x,\frac{n-1}{n})^{-1}h(x,1)],\cr
}
 $$
	and take logarithm of $h(x,\frac{i-1}n)^{-1}h(x,\frac{i}n)$,
	$i=1,\ldots,n$, which is uniformly small if $n$ is large enough.

	(d) Our $h$ is as claimed since we can write $h(x,t)$
	as a so-called product integral, which is a uniform limit as in (b)
	of products of finitely many simple \fns.
	 These products make a multiplicative analog of the ordinary 
	Euler polygon method in linear ordinary differential equations.
	The proof of \p\tEA\ is complete.
}

	 For the rest of this section see [\rPA, \S\S\,7-9], which has all that we
	need but in a slightly different form.
	 We turn now to Runge approximation for simple \fns.
	 Resume the notation and assumptions of \S\,\sC.

\tetel{\p\tEB.}{If $\overline{D_2}$ is \holo{}ally convex in $D_3$,
	$f\in\OOs(\Oa(D_2,R),G)$ is a simple \fn, $0<\ta<1$, and $\epsz>0$,
	then there is a $g\in\OO(D_3\times Y,G)$ that is homotopic through
	\holo\ maps $D_3\times Y\to G$ to\/ $1$ and is \bu\ on any set of the
	form $K\times B_Y(r)$, where\/ $0<r<\infty$ and $K\subsub D_3$, and\/
	$\|f(\za,y)g(\za,y)^{-1}-1\|<\epsz$ for $(\za,y)\in\Oa(D_1,\ta R)$.
}

\biz{Proof.}{Define $f'\in\OOs(\Oa(D_2,R),G)$ and $f''\in\OOs(D_2,G)$
	by writing $f(\za,y)=f'(\za,y)f''(\za)$, where $f''(\za)=f(\za,0)$, so
	$f'(\za,0)=1$.
	 Looking at $f'(\za,ty)$, $t\in[0,1]$, \p\tEA(c) gives
	$f'_i\in\OOs(\Oa(D_2,R),\dot G)$ with
	$f'(\za,y)=\exp(f'_1(\za,y))\ldots\exp(f'_m(\za,y))$.
	 After possibly shrinking $D_2$ arbitrarily slightly
	Grauert's theorem [\rC, Th\'eor\`eme principal (ii)]
	yields $f''_j\in\OOs(D_2,\dot G)$ with
	$f''(\za)=\exp(f''_1(\za))\ldots\exp(f''_n(\za))$.
	 Let $\eta>0$.
	 Apply \p\tCB\ to find \fns\
	$g'_i\in\OO(D_3\times Y,\dot G)$ that approximate $f'_i$ so that
	$\|f'_i(\za,y)-g'_i(\za,y)\|<\eta$ for $(\za,y)\in\Oa(D_1,\ta R)$.
	 After possibly shrinking $D_3$ arbirarily slightly \t{}A of Stein
	theory provides \fns\ $g''_j\in\OOs(D_3,\dot G)$ with
	$\|f''_j(\za)-g''_j(\za)\|<\eta$ for $\za\in D_1$.

	 Letting $g(\za,y)=\exp(g'_1(\za,y))\ldots\exp(g'_m(\za,y))\exp(g''_1(\za))\ldots\exp(g''_n(\za))$
	will do if $\eta>0$ is small enough.
	 The proof of \p\tEB\ is complete.
}

\tetel{\p\tEC.}{Let\/ $8\aa\in\AA$, and choose $\ga\in\AA$ as in \p\tBB(b), and $\ga'\in\AA$
	as in \p\tBB(a).
	Let $f\in\OOs(\Oa_N\<\aa\>,G)$ be a simple \fn, and $\epsz>0$.
\vskip0pt
	{\rm(a)} There is a simple \fn\ $g\in\OO_s(\Oa_{N+1}\<\aa\>,G)$
	with $\|f(x)g(x)^{-1}-1\|<\epsz$ for $x\in\Oa_N\<\ga\>$.
\vskip0pt
	{\rm(b)} There is a $g\in\OO(\Oa,G)$ that is simple on\/
	$\Oa_{N+p}\<\ga'\>$ for all $p\ge0$ with
	$\|f(x)g(x)^{-1}-1\|<\epsz$ for $x\in\Oa_N\<\ga'\>$.
}

\biz{Proof.}{(a) Apply \p\tBB\ and \p\tEB\ as in the proof of
	[\rLC, Thm.\,4.5].

	(b) Letting $g_0=f$, apply (a) repeatedly to find for an $\eta>0$ simple
	\fns\ $g_p\in\OOs(\Oa_{N+p}\<\aa\>,G)$ for $p\ge1$ with
	$\|g_{p-1}(x)g_p(x)^{-1}-1\|<\eta/2^{p+1}$ for $x\in\Oa_{N+p}\<\ga\>$.
	 Setting $g=(\ldots(g_pg_{p-1}^{-1})(g_{p-1}g_{p-2}^{-1})\ldots(g_1g_0^{-1}))g_0=
	\lim_{p\to\infty}g_p$ completes the proof of \p\tEC\ if $\eta>0$ is small enough.
}

	 Now we look at cohomology vanishing for simple \fns.

	 Let $Y$ be a \bspc, $D\subsub D'\subsub D''\subsub\CC^N$ \pscx\ open
	with $\overline{D}$ \holo{}ally convex in $D'$ and $\overline{D'}$
	\holo{}ally convex in $D''$, $R:D''\to(0,\infty)$ \cts\ with $-\log R$ \psh,
	$1/2<\ta<1$, $a<a'<b'<b$ reals, $h\in\OO(D')$,
	$D'_{ab}=D'\cap\{a<\re\,h<b\}$,
	$D_1=D\cap\{\re\,h<b'\}$,
	$D_2=D\cap\{a'<\re\,h\}$,
	$\Oa^\ka=\Oa(D^\ka,\ta R)$ for $\ka=1,2$ as in $(\eCA)$,
	$\Oa^*=\Oa(D'_{ab},R)$,
	$\widetilde\oa\subsub\CC^{N+1}$ a $C^\infty$-smooth strictly smooth
	complete Hartogs domain with
	$\{(\za,\la)\in D\times\CC:|\la|<\ta R(\za)\}
	\subsub\widetilde\oa\subsub
	\{(\za,\la)\in D'\times\CC:|\la|<R(\za)\}$,
	$\widetilde\oa^1=\{(\za,\la)\in\widetilde\oa:\re\,h<b'\}$,
	$\widetilde\oa^2=\{(\za,\la)\in\widetilde\oa:a'<\re\,h\}$,
	$\widetilde\Oa=\{(\za,y)\in\CC^N\times Y:(\za,\|y\|)\in\widetilde\oa\}$,
	$\widetilde\Oa^\ka=\{(\za,y)\in\CC^N\times Y:(\za,\|y\|)\in\widetilde\oa^\ka\}$,
	$\ka=1,2$.
	Let $Z$ be a \bspc, $G=\GL(Z)$.
	Then $\OOs(\widetilde\Oa,Z)$ is a \bspc\ with the sup norm,
	and $\OOs(\widetilde\Oa,G)$ is a \blgp\ with \bla\
	$\OOs(\widetilde\Oa,\dot G)$.

\tetel{\p\tED.}{There are bounded linear operators
	$F_\ka:\OOs(\widetilde\Oa^1\cap\widetilde\Oa^2,Z)\to\OOs(\widetilde\Oa^\ka,Z)$,
	$\ka=1,2$, with $f=F_1(f)+F_2(f)$ for $f\in\OOs(\widetilde\Oa^1\cap\widetilde\Oa^2,Z)$.
}

\biz{Proof.}{See the proof of [\rPA, Prop.\,7.1].
}

	Note that part of \p\tED\ can be reformulated as saying that the
	Mayer--Vietoris sequence
$$
	0\to\OOs(\widetilde\Oa^1\cup\widetilde\Oa^2)
	\buildrel r\over\longrightarrow
	\OOs(\widetilde\Oa^1)\times\OOs(\widetilde\Oa^2)
	\buildrel a\over\longrightarrow
	\OOs(\widetilde\Oa^1\cap\widetilde\Oa^2)\to0
	\tag\eEA
 $$	
	is a split exact sequence of \bspc{}s, where
	$\OOs(\cdot)=\OOs(\cdot,Z)$,
	$r(\Phi)=(\Phi|\widetilde\Oa^1,-\Phi|\widetilde\Oa^2)$, and
	$a(f_1,f_2)=f_1+f_2$.	

\tetel{\p\tEE.}{There are an open \nbd\ ${\Cal N}$ of\/ $1$ in
	$\OOs(\widetilde\Oa^1\cap\widetilde\Oa^2,G)$ and \holo\ maps
	$F_\ka:{\Cal N}\to\OOs(\widetilde\Oa^\ka,G)$, $\ka=1,2$, with
	$F_\ka(1)=1$, and $f=F_1(f)F_2(f)$ for $f\in{\Cal N}$.
}

\biz{Proof.}{Apply \p\tED\ as in the proof of [\rPA, Prop.\,7.2].
}

\tetel{\p\tEF.}{There is an $\epsz_0>0$ \st\ if $f\in\OOs(\Oa^*,G)$ satisfies that
	$\|f(x)-1\|<\epsz_0$ for $x\in\Oa^*$, then there are $f_\ka\in\OOs(\Oa^\ka,G)$,
	$\ka=1,2$, with $f=f_1f_2^{-1}$ on\/ $\Oa^1\cap\Oa^2$.
}

\biz{Proof.}{Apply \p\tEE.
}

\tetel{\p\tEG.}{If $\la\in\OOs(\Oa^*,G)$, $g\in\OOs(\Oa^*,Z)$, then there are
	$g_\ka\in\OOs(\Oa^\ka,Z)$ with $g(x)=g_1(x)-\la(x)g_2(x)$ for $x\in\Oa^1\cap\Oa^2$.
}

\biz{Proof.}{If $\sup_{x\in\Oa^*}\|\la(x)-1\|<\epsz_0$ for an $\epsz_0>0$ small enough,
	then \p\tEF\ gives $\la_\ka\in\OOs(\Oa^\ka,G)$ with $\la=\la_1\la_2^{-1}$.
	We seek $g_\ka$ in the form $g_\ka=\la_\ka h_\ka$, and let $h=\la_1^{-1}g$.
	\p\tDA\ gives $h_\ka\in\OOs(\Oa^\ka,Z)$ with $h=h_1-h_2$, completing the
	proof in the case $\la\approx1$.

	If $\la$ is arbitrary, then we reduce to the case of $\la\approx1$ treated above.
	\p\tEB\ gives a $\La\in\OOs(D''\times B_Y(2\sup_{D'}R),G)$ \st\
	$\|\la(x)\La(x)^{-1}-1\|<\epsz_0$ for $x\in\Oa^*$.
	We seek $g_1,g_2$ in the form $g_1=h_1$, $g_2=\La^{-1}h_2$.
	As the Cousin problem $g(x)=h_1(x)-\la(x)\La(x)^{-1}h_2(x)$ can be
	solved by the already proved case $\la\La^{-1}\approx1$ above, the
	proof of \p\tEG\ is complete.
}

\tetel{\p\tEH.}{If $f\in\OOs(\Oa^*,G)$, then there are $f_\ka\in\OOs(\Oa^\ka,G)$,
	$\ka=1,2$, with $f=f_1f_2^{-1}$.
}

\biz{Proof.}{As in the proof of \p\tEB\ we see that there is a homotopy
	$h\in\OObu(\Oa^*\times B_\CC(2),G)$ with with $h(x,1)=f(x)$
	and $h(x,0)=1$ for $x\in\Oa^*$.
	 We seek to define $f_\ka(x)$ by $f_\ka(x)=h_\ka(x,1)$, where
	$h_\ka\in\OObu(\Oa^\ka\times B_\CC(2\ta),G)$ are to be chosen to satisfy
$$
	h_1(x,t)=h(x,t)h_2(x,t).
\tag\eEB
 $$
 	We will obtain $h_\ka$ from its logarithmic derivative
	$\dot h_\ka=h_{\ka t}h_\ka^{-1}$
	via the parametric initial value problem for ordinary
	differential equations
$$
\left\{
\eqalign{
	\ts\frac{d}{dt}h_\ka(x,t)&=\dot h_\ka(x,t)h_\ka(x,t)\cr
	h_\ka(x,0)&=1\cr
}\right.
.
\tag\eEC
 $$
	Let $\dot h=h_th^{-1}\in\OObu(\Oa^*\times B_\CC(2\ta),\dot G)$ 
	be the left logarithmic derivative of $h$.
	The relation $(\eEB)$ follows by integration of ODEs from $(\eEC)$ together
	with the relation
$$
	\dot h_1=\dot h+h\dot h_2h^{-1},
\tag\eED
 $$
 	which is obtained from $(\eEB)$ by logarithmic differentiation.
	 The equation $(\eED)$ is a Cousin problem
$$
	\dot h=\dot h_1-\la\dot h_2,
\tag\eEE
 $$
	where $\la=\Ad(h)\in\OObu(\Oa^*\times B_\CC(2),\GL(\dot G))$.
	As $(\eEE)$ can be solved by 
	\p\tEG\ the proof of \p\tEH\ is complete.
}

\tetel{\p\tEI.}{If $f\in\OOs^1(\Oa^*,G)$, then there are $f_\ka\in\OOs^1(\Oa^\ka,G)$
	solving $f=f_1f_2^{-1}$ on\/ $(\Oa^1\cap\Oa^2)\times[0,1]$.
}

\biz{Proof.}{Grauert's theorem [\rC, Th\'eor\`eme principal (iii)]
	furnishes $g_\ka\in\OOs^1(D^\ka,G)$ solving
	$f(\za,0,t)=g_1(\za,t)g_2(\za,t)^{-1}$.
	 Defining $f'\in\OOs^1(\Oa^*,G)$ by
	 $f'(\za,y,t)=g_1(\za,t)^{-1}f(\za,y,t)g_2(\za,t)$
	we see that 
	 $f'$ satisfies $f'(\za,0,t)=1$.
	 \p\tEH\ gives $f'_\ka\in\OOs(\Oa^\ka,G)$ with $f'_\ka(\za,0)=1$
	solving $f'(\za,y,1)=f'_1(\za,y)f'_2(\za,y)^{-1}$.
	 Define $f''(\za,y,t)=f_1(\za,ty)^{-1}f'(\za,y,t)f'_2(\za,ty)$.
	 Then $f''(\za,y,t)=1$ for $t\in\{0,1\}$.
	 Let $f''_1(\za,y,t)=f''(\za,\chi(\re\,h(\za))y,t)$,
	$f''_2=1$, where $\chi\in C^\infty(\RR,[0,1])$ is a smooth
	cutoff \fn\ that equals $1$ on $[a',b']$, and $0$ on $\RR\setminus[a,b]$.
	Then setting $f_\ka=g_\ka f'_\ka f''_\ka\in\OOs^1(\Oa^\ka,G)$
	completes the proof of \p\tEI.
}

\alcim{\sF. VANISHING FOR A TRIVIAL BUNDLE.}

	In this section we complete the proof of \t\tAC(a).
	Resume the notation and hypotheses of \t\tAC.

	Let $\UU$ be an open covering of $\Oa$.
	We say that a cochain $f\in C^q(\UU,\OO^Z)$, $q\ge1$, is
	\kiemel{simple} and write $f\in\Cs^q(\UU,\OO^Z)$ if its components
	$f_\sa\in\OOs(|\sa|,Z)$ are simple over the body $|\sa|$ of any $q$-simplex
	$\sa$ of $\UU$.
	The coboundaries and finite sums of simple cochains are simple cochains.
	If $f\in\Cs^q(\UU,\OO^Z)$, and $\da f=0$, then we call $f$ a simple cocycle,
	and write that $f\in\Zs^q(\UU,\OO^Z)$.
	Simple cocycles $\Zs^q(\UU,\OO^Z)$, $q\ge1$, modulo simple coboundaries
	$\da\Cs^{q-1}(\UU,\OO^Z)$ make up a group $\Hs^q(\UU,\OO^Z)$.
	If an open covering $\VV$ of $\Oa$ is any refinement of $\UU$, and $f\in\Cs^q(\UU,\OO^Z)$
	is any simple cochain of $\UU$, then the image of $f$ under any refinement map
	from $\VV$ to $\UU$ is a simple cochain of $\VV$ in $\Cs^q(\VV,\OO^Z)$.
	As $\UU$ gets ever finer it is possible to take the direct limit of the
	$\Hs^q(\UU,\OO^Z)$ to get a limit $\Hs^q(\Oa,\OO^Z)$.
	We claim that $\Hs^q(\Oa,\OO^Z)\cong H^q(\Oa,\OO^Z)$, $q\ge1$, are
	naturally isomorphic.
	To that end it is enough to show that any cochain $f\in C^q(\UU,\OO^Z)$ can be
	represented by a simple cochain $g\in\Cs^q(\VV,\OO^Z)$.
	This can be seen as follows.
	Refine $\UU$ so as to be locally finite ($\Oa$ is a metric space, hence it is
	paracompact).
	Then choose for each point $x\in\Oa$ a ball $V_x=B_X(x,\aa(x))\subset\Oa$ so
	small that if $\sa$ is any of the finitely many $q$-simplices of $\UU$ \st\
	$|\sa|$ intersects $V_x$, then $f_\sa$ and its Fr\'echet differential $df_\sa$
	are bounded on $V_x$.
	Let $\VV=\{V_x:x\in\Oa\}$, and $g$ the image of $f$ under any refinement map
	from $\VV$ to $\UU$.
	Then $\VV$ and $g$ will do the job.
	The most useful case of simple cocycles $\Zs^q(\UU,\OO^Z)$ is that of $q=1$,
	where approximation is necessary, and the extra regularity is welcome.
	If a covering $\UU$ of $\Oa=\Oa(D,R)$ as in $(\eCA)$ has a finite 
	basic refinement $\VV(\Oa)$ as in \S\,\sD, and
	$f\in\Zs^0(\UU,\OO^Z)$ is a simple $0$-cocycle, 
	then $f$ determines a unique simple function
	$f\in\OOs(\Oa,Z)$.
	This is due to the precompactness of $D$.

	 The above has a multiplicative version, too, for $G=\GL(Z)$.
	 Let $\OO_1^G\to\Oa$ be the sheaf of germs of \ho-\fns\ $f:\Oa\times[0,1]\to G$,
	i.e., let $U\subset\Oa$ be open, and define $\OO_1^G(U)=\OO_1(U,G)$
	as the set of all \fns\ $f\in C(U\times[0,1],G)$ with $f(x,0)=1$
	and $f(x,1)$ \holo\ for $x\in U$.
	 We say that a cochain $f\in C^q(\UU,\OO_1^G)$, $q=0,1$, is
	\kiemel{simple} and write that $f\in\Cs^q(\UU,\OO_1^G)$ if its components
	$f_\sa\in\OOs^1(|\sa|,G)$ are simple \ho-\fns\ over the body $|\sa|$ of any
	$q$-simplex $\sa$ of $\UU$.
	 The coboundaries $f_Uf_V^{-1}$, $f_{UV}f_{VW}f_{WU}$, and finite products
	and inverses of simple cochains are simple cochains.
	 If $f=(f_{UV})\in\Cs^1(\UU,\OO_1^G)$ and $\da f=(f_{UV}f_{VW}f_{WU})=1$,
	then we call $f$ a simple cocycle, and write that $f\in\Zs^1(\UU,\OO_1^G)$.
	 Between two simple cochains $f=(f_{UV}),g=(g_{UV})\in\Cs^1(\UU,\OO^G_1)$
	there is an equivalence relation $f\sim g$ defined by making $f\sim g$
	if and only if there is a simple cochain $c\in\Cs^0(\UU,\OO_1^G)$ with
	$f_{UV}=c_U^{-1}g_{UV}c_V$.
	 Simple cocycles $\Zs^1(\UU,\OO^G_1)$ modulo this equivalence
	relation $\sim$ make up a simple cohomology set $\Hs^1(\UU,\OO^G_1)$
	with a distiguished element $1$, which is just the class of the trivial
	cocycle $1$.
	 If an open covering $\VV$ of $\Oa$ is any refinement of $\UU$, and
	$f\in\Cs^q(\UU,\OO^G_1)$ is any simple cochain, then the image of $f$
	under any refinement map $\VV\to\UU$ is a simple cochain of $\VV$ in
	$\Cs^q(\VV,\OO^G_1)$, $q=0,1$.
	 As $\UU$ gets ever finer it is possible to take the direct limit of the
	$\Hs^1(\UU,\OO^G_1)$ to get the simple cohomology group $\Hs^1(\Oa,\OO^G_1)$
	of $\Oa$ \wrt\ the sheaf $\OO^G_1$.
	If a covering $\UU$ of $\Oa=\Oa(D,R)$ as in $(\eCA)$ has a finite 
	basic refinement $\VV(\Oa)$ as in \S\,\sD, and
	$f\in\Zs^0(\UU,\OO^G_1)$ is a simple $0$-cocycle, 
	then $f$ determines a unique simple \ho-\fn\
	$f\in\OOs^1(\Oa,G)$.
	This is due to the precompactness of $D$.

\tetel{\p\tFA.}{Let $X$ be a \bspc\ with a \sbs, $\Oa\subset X$ \pscx\ open,
	and suppose that \pshdom\ holds in $\Oa$.
	Then for any $\aa\in\AA$ there is a $\ga\in\AA$ \st\ $\ga<\aa$, and
	$\Hs^q(\BB_N(\aa),\OO^Z)|\BB_N(\ga)=0$ for all $N\ge0$ and $q\ge1$.
}

\biz{Proof.}{We consider some open coverings and refinement maps of them.
	Let $\aa,\ba,\ga\in\AA$ be as in \p\tBB(b).
	Consider the open coverings
	$\BB_N(\aa)$,
	$\UU_N=\{U(x)=\Oa_N\<\ba\>\cap\pi_N^{-1}B_X(\pi_Nx,\ba(x)):x\in\Oa_N\<\ga\>\}$,
	$\BB_N(\ga)$, and their refinement maps
	$\UU_N\to\BB_N(\aa)$ given by $U(x)\mapsto B_X(x,\aa(x))$, and
	$\BB_N(\ga)\to\UU_N$ given by $B_X(x,\ga(x))\mapsto U(x)$.
	Due to the inequalities $(\eBB)$ the above are indeed refinement maps, and
	hence induce maps
$$
	\Hs^q(\BB_N(\aa),\OO^Z)\to\Hs^q(\UU_N,\OO^Z)\to\Hs^q(\BB_N(\ga),\OO^Z)
 $$
	in cohomology for $N\ge0$ and $q\ge1$.
	Since the first map has zero image by \p\tDA\ the proof of \p\tFA\ is
	complete.
}

\biz{Proof of \t\tAC(a).}{Let $f\in H^q(\Oa,\OO^Z)=\Hs^q(\Oa,\OO^Z)$ be a
	cohomology class that we would like to resolve.
	By \pshdom\ in $\Oa$ there is an $\aa$ \st\ $10\aa\in\AA$, and $f$ can be
	represented by a simple cocycle $f\in\Zs^q(\BB(\aa),\OO^Z)$.
	On choosing a $\ga\in\AA$ as in \p\tFA\ we find
	$g_N\in\Cs^{q-1}(\BB_N(\ga),\OO^Z)$, $N\ge0$, with $\da g_N=f|\BB_N(\ga)$.
	We can extend the cochain $g_N$ to a cochain $g_N\in\Cs^{q-1}(\BB(\ga),\OO^Z)$
	simply by defining $g_N$ to be zero over simplices $\bigcap_{i=1}^q B_X(x_i,\ga(x_i))$
	if at least one vertex $x_i\not\in\Oa_N\<\ga\>$.
	\p\tBB(e) gives a $\ba\in\AA$, $\ba<\ga$, with 
	$\Oa_N\<\ba\>\subset\Oa_N\<\ga\>\cap\Oa_{N+1}\<\ga\>$ for $N\ge0$.
	So $(g_{N+1}-g_N)|\BB_N(\ba)\in\Zs^q(\BB_N(\ba),\OO^Z)$.

	Suppose first that $q\ge2$.
	\p\tFA\ gives a $\ba'\in\AA$ \st\ $\ba'<\ba$, and
	$\Hs^{q-2}(\BB_N(\ba),\OO^Z)|\BB_N(\ba')=0$,
	so similarly extending a $(q-2)$-cochain there is an
	$h_N\in\Cs^{q-2}(\BB(\ba'),\OO^Z)$ with
	$(g_{N+1}-g_N)|\BB_N(\ba')=\da h_N|\BB_N(\ba')$.
	Letting $g'_N=g_N|\BB(\ba')-\sum_{n=1}^{N-1}\da h_n\in\Cs^{q-1}(\BB(\ba'),\OO^Z)$
	\p\tBB(a) implies as $g'_{N+1}|\BB_N(\ba')=g_N|\BB_N(\ba')$ that
	$g'_N$ converges as $N\to\infty$ to a cochain
	$g\in C^{q-1}(\BB(\ba'),\OO^Z)$ with $\da g=f|\BB(\ba')$.
	Thus $f$ equals zero in $H^q(\Oa,\OO^Z)$ for $q\ge2$.

	Let now $q=1$.
	By \p\tBB(f) there is a $\ba'\in\AA$, $\ba'<\ba$, \st\ the covering
	$\BB_N(\ba)|\Oa_N\<\ba'\>$ has a finite refinement for all $N\ge0$.
	As $(g_{N+1}-g_N)|\BB_N(\ba)=h_N\in\Zs^0(\BB_N(\ba),\OO^Z)$
	we see that over $\Oa_N\<\ba'\>$ our $h_N|(\BB_N(\ba)|\Oa_N\<\ba'\>)$
	patches up to simple \fn\ $h_N\in\OOs(\Oa_N\<\ba'\>,Z)$.
	\p\tCC(b) gives a $\ba''\in\AA$, $\ba''<\ba'$, and an
	$H_N\in\OO(\Oa,Z)$ with
	$\|h_N(x)-H_N(x)\|<1/2^N$ for $x\in\Oa_N\<\ba''\>$ for all $N\ge0$.
	\p\tBB(a) yields a $\ba'''\in\AA$, $\ba'''<\ba''$ \st\
	$\Oa_M\<\ba'''\>\subset\Oa_N\<\ba''\>$ for all $N\ge M\ge0$.
	Letting $g'_N=(g_N-\sum_{n=1}^{N-1}H_n)|\BB(\ba'')$
	defines a cochain $g'_N\in C^0(\BB(\ba''),\OO^Z)$ that satisfies
	$g_1+\sum_{n=1}^\infty(g_{n+1}-g_n-H_n)=
	\lim_{N\to\infty}(g_N-\sum_{n=1}^{N-1}H_n)=\lim_{N\to\infty}g'_N=g$,
	where the convergence is uniform on $\Oa_M\<\ba'''\>$ for all $M\ge0$.
	Thus the limit $g\in C^0(\BB(\ba''),\OO^Z)$ exists and satisfies
	$\da g=\da g_N|\BB(\ba'')=f|\BB(\ba'')$.
	The proof of \t\tAC(a) is complete.
}

	 Resume the context and notation of \S\,\sD.

\tetel{\p\tFB.}{For any $f\in\Zs^1(\UU'(\Oa'),\OO^G)$ with $f(\za,0)=1$
	there is a $g\in\Cs^0(\UU(\Oa),\OO^G)$ with
	$g(\za,0)=1$ \st\
	$f_{U'V'}(\za,y)|(U\cap V)(\Oa)=g_U(\za,y)g_V(\za,y)^{-1}$.
}

\biz{Proof.}{Letting $h\in\OO(\CC^N)$ be various linear \fns\ $h(\za_1,\ldots,\za_N)=\za_j$
	the usual induction process of Cousin and Cartan (see [\rH, \S\,7.2],
	or [\rLt, Lemma\,4.1]) relying on \p\tEH\ completes the proof of \p\tFB.
}

\tetel{\p\tFC.}{For any $f\in\Zs^1(\UU'(\Oa'),\OO^G_1)$ with $f(\za,0,t)=1$
	there is a $g\in\Cs^0(\UU(\Oa),\OO^G_1)$ with $g(\za,0,t)=1$ \st\
	$f_{U'V'}(\za,y,t)|(U\cap V)(\Oa)=g_U(\za,y,t)g_V(\za,y,t)^{-1}$.
}

\biz{Proof.}{The proof is a similar induction as that of \p\tFB\ relying on
	\p\tEI\ this time.
}

\tetel{\p\tFD.}{For any $f\in\Zs^1(\UU'(\Oa'),\OO^G_1)$ there is a
	$g\in\Cs^0(\UU(\Oa),\OO^G_1)$ with
	$f_{U'V'}(\za,y,t)|(U\cap V)(\Oa)=g_U(\za,y,t)g_V(\za,y,t)^{-1}$.
}

\biz{Proof.}{Define $f'\in\Zs^1(\UU',\OO^G_1)$ and $f''\in\Zs^1(\UU'(\Oa'),\OO^G_1)$
	by writing $f(\za,y,t)=f'(\za,t)f''(\za,y,t)$, where $f'(\za,t)=f(\za,0,t)$,
	so $f''(\za,0,t)=1$.
	 After possibly shrinking $\UU'$ arbitrarily slightly
	 Grauert's theorem [\rC, Th\'eor\`eme principal (iii)]
	gives a $g'\in\Cs^0(\UU',\OO^G_1)$ with
	$f'_{U'V'}(\za,t)=g'_{U'}(\za,t)g'_{V'}(\za,t)^{-1}$.
	 Define a new cocycle $f''\in\Zs^1(\UU'(\Oa'),\OO^G_1)$ by the formula
	$f''_{U'V'}(\za,y,t)=g'_{U'}(\za,t)^{-1}f_{U'V'}(\za,y,t)g'_{V'}(\za,t)$.
	 Then $f''(\za,0,t)=1$, and \p\tFC\ provides a
	$g''\in\Cs^0(\UU(\Oa),\OO^G_1)$ with
	$g''(\za,0,t)=1$, and
	$f''_{U'V'}(\za,y,t)|(U\cap V)(\Oa)=g''_U(\za,y,t)g''_V(\za,y,t)^{-1}$.
	 Letting $g_{U}(\za,y,t)=g'_{U'}(\za,t)g''_U(\za,y,t)$
	completes the proof of \p\tFD.
}

	 We now turn to approximation.
	 Resume the context and notation of \S\,\sC\ and \p\tEB.

\tetel{\p\tFE.}{Let $\overline{D_2}$ be \holo{}ally convex in $D_3$, $0<\ta<1$,
	$\epsz>0$, and $f\in\OOs^1(\Oa(D_2,R),G)$.
\vskip0pt
	{\rm(a)} If $f(\za,0,t)=1$, then there is a $g\in\OO_1(D_3\times Y,G)$
	with $g(\za,0,t)=1$ \st\ $g(\za,y,t)$ is \bu\ on any set of the form
	$K\times B_Y(r)\times[0,1]$, where\/ $0<r<\infty$, and $K\subsub D_3$,
	and\/ $\|f(\za,y,t)g(\za,y,t)^{-1}-1\|<\epsz$ for $(\za,y)\in\Oa(D_1,\ta R)$,
	and $t\in[0,1]$.
\vskip0pt
	{\rm(b)} There is a $g\in\OO_1(D_3\times Y,G)$ \st\ $g(\za,y,t)$
	is \bu\ on any set of the form
	$K\times B_Y(r)\times[0,1]$, where\/ $0<r<\infty$, and $K\subsub D_3$,
	and\/ $\|f(\za,y,t)g(\za,y,t)^{-1}-1\|<\epsz$ for $(\za,y)\in\Oa(D_1,\ta R)$,
	and $t\in[0,1]$.
}

\biz{Proof.}{(a) Looking back at the proof of \p\tEB\ we find that it gives
	a $g'\in\OO(D_3\times Y,G)$ with $g'(\za,0)=1$ and
	$\|f(\za,y,1)g'(\za,y)^{-1}-1\|<\epsz$ for $(\za,y)\in\Oa(D_1,\ta R)$.
	Let $\chi\in C^\infty(\CC^N,[0,1])$ be a smooth cutoff \fn\
	with $\chi=1$ on $D_1$, and $\chi=0$ on $\CC^N\setminus D_2$.
	Fix $\ta'$ with $0<\ta<\ta'<1$, and define
	$\tilde f\in C(D_3\times Y\times[0,1],G)$ by
	$\tilde f(\za,y,t)=f(\za,\chi(\za)\min(\|y\|,\ta'R(\za))\frac{y}{\|y\|},t)$.
	Let $g(\za,y,t)=g'(\za,ty)\tilde f(\za,ty,1)^{-1}\tilde f(\za,y,t)$.
	We check that this $g$ will do.
	Indeed, $g$ is defined and \cts\ on $D_3\times Y\times[0,1]$,
	$g(\za,y,0)=1$, $g(\za,0,t)=1$, $g(\za,y,1)=g'(\za,y)$ is \holo\
	on $D_3\times Y$, and for $(\za,y)\in\Oa(D_1,\ta R)$, $t\in[0,1]$ we have that
$
	\|f(\za,y,t)g(\za,y,t)^{-1}-1\|\le\|f(\za,ty,1)g'(\za,ty)^{-1}-1\|<\epsz
 $
 	as $f=\tilde f$ on $\Oa(D_1,\ta R)\times[0,1]$, and $(\za,ty)\in\Oa(D_1,\ta R)$.

	(b) Let $\eta>0$, and write $f'(\za,t)=f(\za,0,t)$, $f(\za,y,t)=f'(\za,t)f''(\za,y,t)$;
	$f'\in\OOs^1(D_2,G)$, $f''\in\OOs^1(\Oa(D_2,R),G)$.
	 Grauert's theorem [\rC, Th\'eor\`eme principal (ii)] gives a $g'\in\OOs^1(D_3,G)$
	with $\|f'(\za,t)g'(\za,t)^{-1}-1\|<\eta$ for $(\za,t)\in D_1\times[0,1]$.
	 Part (a) supplies a $g''\in\OO_1(D_3\times Y,G)$ with
	$g''(\za,0,t)=1$ and $\|f''(\za,y,t)g''(\za,y,t)^{-1}-1\|<\eta$
	for $(\za,y)\in\Oa(D_1,\ta R)$, and $t\in[0,1]$.
	 Letting $g=g'g''$ completes the proof of \p\tFE\ if $\eta>0$ is small enough.
}

\tetel{\p\tFF.}{Let\/ $8\aa\in\AA$, and choose a $\ga\in\AA$ as in \p\tBB(b), and
	a $\ga'\in\AA$ as in \p\tBB(a).
	 Let $f\in\OOs^1(\Oa_N\<\aa\>,G)$ be a simple \ho-\fn, and $\epsz>0$.
\vskip0pt
	{\rm(a)} There is a simple \ho-\fn\ $g\in\OOs^1(\Oa_{N+1}\<\aa\>,G)$ 
	that satisfies
	$\|f(x,t)g(x,t)^{-1}-1\|<\epsz$ for $x\in\Oa_N\<\ga\>$, $t\in[0,1]$.
\vskip0pt
	{\rm(b)} There is a $g\in\OO_1(\Oa,G)$ that is a simple \ho-\fn\
	on\/ $\Oa_{N+p}\<\ga'\>$ for all $p\ge0$ with
	$\|f(x,t)g(x,t)^{-1}-1\|<\epsz$ for $x\in\Oa_N\<\ga'\>$, $t\in[0,1]$.
}

\biz{Proof.}{Relying on \p\tFE(b) instead of \p\tEB\ as in the proof of \p\tEC\
	completes the proof of \p\tFF.
}

\tetel{\p\tFG.}{Let $X$ be a \bspc\ with a \sbs, $\Oa\subset X$ \pscx\ open,
	and suppose that \pshdom\ holds in $\Oa$.
	Then for any $\aa\in\AA$ there is a $\ga\in\AA$ \st\ $\ga<\aa$, and
	$\Hs^1(\BB_N(\aa),\OO^G_1)|\BB_N(\ga)=1$ for all $N\ge0$.
}

\biz{Proof.}{We consider some open coverings and refinement maps of them.
	Let $\aa,\ba,\ga\in\AA$ be as in \p\tBB(b).
	Consider the open coverings
	$\BB_N(\aa)$,
	$\UU_N=\{U(x)=\Oa_N\<\ba\>\cap\pi_N^{-1}B_X(\pi_Nx,\ba(x)):x\in\Oa_N\<\ga\>\}$,
	$\BB_N(\ga)$, and their refinement maps
	$\UU_N\to\BB_N(\aa)$ given by $U(x)\mapsto B_X(x,\aa(x))$, and
	$\BB_N(\ga)\to\UU_N$ given by $B_X(x,\ga(x))\mapsto U(x)$.
	Due to the inequalities $(\eBB)$ the above are indeed refinement maps, and
	hence induce maps
$$
	\Hs^1(\BB_N(\aa),\OO^G_1)\to\Hs^1(\UU_N,\OO^G_1)\to\Hs^1(\BB_N(\ga),\OO^G_1)
 $$
	in cohomology for $N\ge0$.
	Since the image of the first map is $1$ by \p\tFD\ the proof of \p\tFG\ is
	complete.
}

\tetel{\t\tFH.}{Let $X$ be a \bspc\ with a \sbs, $\Oa\subset X$ \pscx\ open,
	and suppose that \pshdom\ holds in $\Oa$.
	Then $\Hs^1(\Oa,\OO^G_1)=1$.
}

\biz{Proof.}{Let $f\in\Hs^1(\Oa,\OO^G_1)$ be a cohomology class that we would like to
	resolve.
	 By \pshdom\ in $\Oa$ there is an $\aa$ \st\ $10\aa\in\AA$, and
	$f$ can be represented by a simple cocycle $f\in\Zs^1(\BB(\aa),\OO^G_1)$.
	 On choosing a $\ga\in\AA$ as in \p\tFG\ we find $g^N\in\Cs^0(\BB_N(\ga),\OO^G_1)$,
	$N\ge0$, with $f_{UV}|\BB_N(\ga)=g^N_U(g^N_V)^{-1}$.
	 We can extend the cochain $g^N$ to a cochain $g^N\in\Cs^0(\BB(\ga),\OO^G_1)$
	simply by defining $g^N_U$ to be $1$ over simplices $U=B_X(x,\ga(x))$
	if $x\not\in\Oa_N\<\ga\>$.
	 \p\tBB(e) gives a $\ba\in\AA$, $\ba<\ga$, with $\Oa_N\<\ba\>\subset\Oa_N\<\ga\>\cap
	\Oa_{N+1}\<\ga\>$ for $N\ge0$.
	 So $((g^N)^{-1}g^{N+1})|\BB_N(\ba)\in\Zs^0(\BB(\ba),\OO^G_1)$.
	 By \p\tBB(f) there is a $\ba'\in\AA$, $\ba'<\ba$, \st\
	the covering $\BB_N(\ba)|\Oa_N\<\ba'\>$ has a finite basic refinement for all $N\ge0$.
	 As $((g^N)^{-1}g^{N+1})|\BB_N(\ba)=h_N\in\Zs^0(\BB_N(\ba),\OO^G_1)$ we see that
	over $\Oa_N\<\ba'\>$ our $h_N|(\BB_N(\ba)|\Oa_N\<\ba'\>)$ patches up to a simple
	\ho-\fn\ $h_N\in\OOs^1(\Oa_N\<\ba'\>,G)$.
	 Let $a_1=1$.
	 Repeated application of \p\tFF(b) gives a $\ba''\in\AA$, $\ba''<\ba'$, and
	a sequence $\tilde g^N\in C^0(\BB(\ba''),\OO^G_1)$, $N\ge1$, of the form
	$\tilde g^N_U=g^N_Ua_N$, where $a_N\in\OO_1(\Oa,G)$ is \st\ 
	$a_N|\Oa_{N+p}\<\ba''\>\in\OOs^1(\Oa_{N+p}\<\ba''\>,G)$ is simple for all
	$p\ge0$, and
	$(\tilde g^N)^{-1}\tilde g^{N+1}=a_N^{-1}h_Na_{N+1}$ satisfies that
	$\|a_N^{-1}h_Na_{N+1}-1\|<1/2^N$ on $\Oa_N\<\ba'''\>$.
	 \p\tBB(a) yields a $\ba'''\in\AA$, $\ba'''<\ba''$ with
	$\Oa_M\<\ba'''\>\subset\Oa_N\<\ba''\>$ for all $N\ge M\ge0$.
	 As $\tilde g^N$ converges uniformly on $\Oa_M\<\ba'''\>$ for all $M$
	the limit $g=\lim_{N\to\infty}\tilde g^N$ exists and satisfies that
	$g\in C^0(\BB(\ba''),\OO^G_1)$, $g|\BB(\ba''')\in\Cs^0(\BB(\ba'''),\OO^G_1)$,
	and $g_Ug_V^{-1}=f_{UV}|\BB(\ba'')$.
	 Thus $\Hs^1(\Oa,\OO^G_1)=1$, and the proof of \t\tFH\ is complete.
}

\alcim{\sG. TOPOLOGICAL TRIVIALITY OF A VECTOR BUNDLE.}

	 In this section we show that certain type of \tpl\ \bvbdl{}s are
	\tpl{}ly trivial.

	 Given a \bspc\ $Z$, and a number $1\le p<\infty$, let
	$Z_1=\ell_p(Z)=\{z=(z_n):z_n\in Z,\|z\|=(\sum_{n=1}^\infty \|z_n\|^p)^{1/p}<\infty\}$.
	Note that $Z_1$, $Z\oplus Z_1$, $Z_1\oplus Z_1$, $\ell_p(Z_1)$ are
	isomorphic \bspc{}s by permuting the coordinates.
	 Denote by $Z_1$ the trivial \bvbdl\ $\Oa\times Z_1$ over any base space $\Oa$.

\tetel{\p\tGA.}{Let $X,Z$ be \bspc{}s, $X$ with a \sbs, $\Oa\subset X$ open,
	$E\to\Oa$ a \tpl\ \bvbdl\ with fiber type $Z$, and $F=E\oplus Z_1\to\Oa$
	the direct sum \bvbdl.
	 Then $F$ is \tpl{}ly trivial: $F$ is \cts{}ly isomorphic to\/ $\Oa\times Z_1$.
}

	 \p\tGA\ is the main point of this section, and its proof will take us
	some steps.

	 Let $\NN=\{1,2,3,\ldots\}$, and fix bijections from $\NN$ to members of
	a partition of $\NN$ into two infinite subsets, say, odd or even numbers.
	 Let $I_2$ be the isomorphism $I_2\in\Hom(Z_1,Z_1\oplus Z_1)$ of \bspc{}s
	induced by the above bijection, i.e., $I_2((z_n))=((x_n),(y_n))$,
	where $x_n=z_{2n-1}$, and $y_n=z_{2n}$ for $n\ge1$.
	 Fix a bijection $\NN\to\NN\times\NN$, say, $n\mapsto(i,j)$, where
	$n=2^{i-1}(2j-1)$.
	 Let $I_\infty$ be the isomorphism $I_\infty\in\Hom(Z_1,\ell_p(Z_1))$
	of \bspc{}s induced by the above bijection, i.e.,
	$I_\infty((z_n))=(i\mapsto(j\mapsto x_{ij}))$, where $x_{ij}=z_n$, and
	$n=2^{i-1}(2j-1)$.
	 Then $I_2,I_\infty$ are isometries and their operator norms are
	$\|I_2\|=\|I_\infty\|=1$.

	 Our \vbdl\ $F$ has the stability property that $F$ and 
	$F\oplus Z_1=E\oplus Z_1\oplus Z_1$ are
	isomorphic.
	 This is due to the fact that $Z_1$ and $Z_1\oplus Z_1$ are isomorphic.
	 Let $J\in C(\Oa,\Hom(F,F\oplus Z_1))$ be an isomorphism of \tpl\ \bvbdl{}s
	defined, e.g., by $J(\xi,z_1)=(\xi,I_2(z_1))$, where $\xi\in E_x$, $x\in\Oa$,
	and $z_1\in Z_1$.

	 Denote a finite or infinite block diagonal matrix $A$ with
	diagonal blocks $A_1,A_2,A_3,\ldots$ by
	$A=\diag(A_1,A_2,A_3,\ldots)$.

\tetel{\p\tGB.}{Let $U\subset\Oa$ be open, $f\in C(U,\GL(Z_1))$, and
	define $f'\in C(U,\GL(Z_1\oplus Z_1))$ by
	$f'(x)=\diag(f(x),1)$.
	 Then $f'$ is null homotopic, i.e., there is an
	$h\in C(U\times[0,1],\GL(Z_1\oplus Z_1))$ with
	$h(x,0)=f'(x)$, and $h(x,1)=1$ for all $x\in U$.
}

\biz{Proof.}{This is based on classical tricks with infinite matrices,
	see [\rK, Lemma\,7].
	 As $Z_1\oplus Z_1$, $Z_1\oplus\ell_p(Z_1)$, and $\ell_p(Z_1)$
	are isomorphic by 
	$Z_1\oplus Z_1\ni(x,y)\mapsto(x,I_\infty(y))\in Z_1\oplus\ell_p(Z_1)$,
	and by $Z_1\oplus\ell_p(Z_1)\ni(x,(y_n))\mapsto(z_n)\in\ell_p(Z_1)$,
	where $z_1=x$, and $z_n=y_{n-1}$ for $n\ge2$, we see that $f'$ can be
	regarded as an element of $C(U,\GL(\ell_p(Z_1)))$ defined by the
	infinite block diagonal matrix $f'(x)=\diag(f(x),1,1,\ldots)$.
	
	 Our homotopy $h$ will be the concatenation $h=h_1\lor h_2$ of two
	homotopies $h_i\in C(U\times[0,1],\GL(\ell_p(Z_1)))$, $i=1,2$, defined
	by $h(x,t)=h_1(x,2t)$ for $x\in U$, $t\in[0,\frac12]$, and 
	$h(x,t)=h_2(x,2t-1)$ for $x\in U$, $t\in[\frac12,1]$, where we must have
	$h_1(x,1)=h_2(x,0)$ for all $x\in U$.

	 To define $h_1$ let $h_1(x,s)$ equal
$$
	\diag\Big(f(x),
	\bmatrix
	\cos t&-\sin t\cr
	\sin t&\cos t\cr
	\endbmatrix
	\bmatrix
	f(x)&0\cr
	0&1\cr
	\endbmatrix
	\bmatrix
	\cos t&\sin t\cr
	-\sin t&\cos t\cr
	\endbmatrix
	\bmatrix
	f(x)^{-1}&0\cr
	0&1\cr
	\endbmatrix,
	\ldots\Big),
 $$
 	where the second diagonal block repeats along the diagonal, and $t=\frac\pi2s$.
	 Then we can easily compute that $h_1(x,0)=\diag(f(x),1,1,1,\ldots)$, and
	$h_1(x,1)=\diag(f(x),f(x)^{-1},f(x),f(x)^{-1},\ldots)$.

	 We define $h_2$ by letting $h_2(x,s)$ equal
$$
        \diag\Big(
        \bmatrix
        \cos t&-\sin t\cr
        \sin t&\cos t\cr
        \endbmatrix
        \bmatrix
        f(x)^{-1}&0\cr
        0&1\cr
        \endbmatrix
        \bmatrix
        \cos t&\sin t\cr
        -\sin t&\cos t\cr
        \endbmatrix
        \bmatrix
        f(x)&0\cr
        0&1\cr
        \endbmatrix,
        \ldots\Big),
 $$
	where the first diagonal block repeats along the diagonal, and
	$t=\frac\pi2(s+1)$.
	 Then $h_2(x,0)=h_1(x,1)$, and $h_2(x,1)=\diag(1,1,1,1,\ldots)$.
	 The proof of \p\tGB\ is complete.
}

	 Elementary algebraic topology tells us that null homotopic maps have
	an extension property.

\tetel{\p\tGC.}{For any $F\subset U\subset\Oa$, with $F$ closed, $U$ open, and
	$f\in C(U,\GL(Z_1))$ there is an $f''\in C(\Oa,\GL(Z_1\oplus Z_1))$
	with $f''(x)=f'(x)$ for $x\in F$, where $f'$ is as in \p\tGB.
}

\biz{Proof.}{\p\tGB\ gives a homotopy $h\in C(U\times[0,1],\GL(Z_1\oplus Z_1))$
	with $h(x,0)=f'(x)$ and $h(x,1)=1$ for all $x\in U$.
	 Let $\chi\in C(\Oa,[0,1])$ be a cutoff \fn\ that equals $1$ on $F$ and
	$0$ on $\Oa\setminus U$.
	 Then letting $f''(x)$ equal $h(x,1-\chi(x))$ for $x\in U$ and $1$ for
	$x\in\Oa\setminus U$ completes the proof of \p\tGC.
}

\tetel{\p\tGD.}{Let $U_1,U_2\subset\Oa$ be open, and suppose that there are
	trivializations $t_i\in C(U_i,\Hom(F,Z_1))$ of $F|U_i$, $i=1,2$.
	 Then $F\oplus Z_1$ and $F$ are trivial over $U_1\cup U_2$.
}

\biz{Proof.}{Let $U=U_1\cap U_2$, and define $f\in C(U,\GL(Z_1))$ by
	$f(x)=t_2(x)t_1(x)^{-1}$, and $f'\in C(U,\GL(Z_1\oplus Z_1))$
	by $f'(x)=\diag(f(x),1)$.
	 Let $V_i\subset\overline{V_i}\subset U_i$ be a shrinking of the
	covering $\{U_1,U_2\}$ relative to $U_1\cup U_2$.
	 \p\tGC\ gives an extension $f''\in C(\Oa,\GL(Z_1\oplus Z_1))$
	with $f''(x)=f'(x)$ for $x\in\overline{V_1}\cap\overline{V_2}$.
	 Define the trivializations $T_i\in C(V_i,\Hom(F\oplus Z_1,Z_1\oplus Z_1))$
	by $T_1(x)=f''(x)\diag(t_1(x),1)$, and $T_2(x)=\diag(t_2(x),1)$.
	As $T_1(x)=T_2(x)$ for $x\in V_1\cap V_2$, and $V_1\cup V_2=U_1\cup U_2$,
	the $T_1$ and $T_2$ patch up to define a trivialization $T\in C(U_1\cup U_2,
	\Hom(F\oplus Z_1,Z_1\oplus Z_1))$ by letting $T(x)$ equal $T_i(x)$ for $x\in V_i$,
	$i=1,2$.
	 Hence we get that $F\oplus Z_1$ is trivial over $U_1\cup U_2$.
	 Defining $S\in C(U_1\cup U_2,\Hom(F,Z_1))$ by
	$S(x)=I_2^{-1}T(x)J(x)$ trivializes $F|U_1\cup U_2$ as well.
	 The proof of \p\tGD\ is complete.
}

	 Thus $F$ is trivial over any finite union $\bigcup U_i$ of open subsets
	$U_i$ of $\Oa$ over each of which $F|U_i$ is trivial.

\biz{Proof of \p\tGA.}{Let $\AA'$ as in $\S\,\sB$.
	There is an $\aa\in\AA'$ \st\ $F|B_X(x,\aa(x))$ is trivial for all $x\in\Oa$.
	\p\tBB\ (except its part (c)) remains true for any $\aa\in\AA'$ as well.
	The proof is just like that of \p\tBB\ for $\AA$ only it is easier,
	and \pshdom\ in $\Oa$ is not needed.
	\p\tBB(f) gives an $\aa'\in\AA'$, $\aa'<\aa$, \st\ $\BB_N(\aa)|\Oa_N(\aa')$
	has a finite refinement for all $N\ge0$.
	\p\tBB(e) yields a $\ba\in\AA'$, $\ba<\aa'$, with $\Oa_N\<\ba\>\subset
	\Oa_N\<\aa'\>\cap\Oa_{N+1}\<\aa'\>$ for all $N\ge0$, and \p\tBB(a) a
	$\ga\in\AA'$, $\ga<\ba$, with $\overline{\Oa_n\<\ga\>}\subset\Oa_N\<\ba\>$
	for $N\ge n$.

	 We show by induction that there are for any $N\ge1$ \cts\ trivializations
	$T_N\in C(\Oa_N\<\aa'\>,\Hom(F\oplus\ell_p(Z_1),\ell_p(Z_1)))$
	of the form $T_N(x)=\diag(t_N(x),1,1,\ldots)$, where
	$t_N\in C(\Oa_N\<\aa'\>,\Hom(F\oplus Z_1^{\oplus 2^{N-1}-1}, Z_1^{\oplus 2^{N-1}}))$
	are trivializations \st\ $T_N=T_{N+1}$ over $\Oa_M\<\ga\>$ for all $N\ge M\ge0$.

	 By the remark after the proof of \p\tGD\ there is a trivialization
	$t_1\in C(\Oa_1\<\aa'\>,\Hom(F,Z_1))$.
	 Define $T_1$ by $T_1(x)=\diag(t_1(x),1,1,\ldots)$.
	 Suppose that $T_1,\ldots, T_N$ are already defined.
	 By the remark after the proof of \p\tGD\ there is a trivialization
	$T'_{N+1}\in C(\Oa_{N+1}\<\aa'\>,\Hom(F\oplus\ell_p(Z_1),\ell_p(Z_1)))$
	of the form $T'_{N+1}(x)=\diag(t'_{N+1}(x),1,1,\ldots)$, where
	$t'_{N+1}\in C(\Oa_{N+1}\<\aa'\>,\Hom(F\oplus Z_1^{\oplus2^N-1},Z_1^{\oplus2^N}))$
	is a trivialization.
	 We define $f\in C(\Oa_N\<\aa'\>\cap\Oa_{N+1}\<\aa'\>,\GL(Z_1^{\oplus2^N}))$
	by $f(x)=t_N(x)t'_{N+1}(x)^{-1}$, and 
	$f'\!\in\! C(\Oa_N\<\aa'\>\cap\Oa_{N+1}\<\aa'\>,\GL(Z_1^{\oplus2^N}\!\!\oplus Z_1^{\oplus2^N}))$
	by $f'(x)=\diag(f(x),1)$.
	 Then by \p\tGB\ our $f'$ is null homotopic, and thus by \p\tGC\
	there is an extension $f''\in C(\Oa,\GL(Z_1^{\oplus2^N}\oplus Z_1^{\oplus2^N}))$
	of $f'$ with $f''(x)=f'(x)$ for $x\in\overline{\Oa_N\<\ba\>}$.

	 Define $T_{N+1}$ by $T_{N+1}(x)=\diag(f''(x),1,1,\ldots)T'_{N+1}$.
	 As $T_{N+1}(x)=\diag(f'(x),1,1,\ldots)\diag(t'_{N+1},1,1,\ldots)=
	 (t_N(x),1,1,\ldots)=T_N(x)$ for all $x\in\Oa_M\<\ga\>$ for $N\ge M\ge0$
	 we see that the sequence $T_N$ converges as $N\to\infty$ to a trivialization
	 $T\in C(\Oa,\Hom(F\oplus\ell_p(Z_1),\ell_p(Z_1)))$ in a quasi-stationary
	 manner over $\Oa_M\<\ga\>$ for all $M$.
	  As $F\oplus\ell_p(Z_1)\cong F$, and $\ell_p(Z_1)\cong Z_1$
	 the proof of \p\tGA\ is complete.
	
}

\alcim{\sH. FROM CONTINUOUS TO LOCALLY LIPSCHITZ SECTIONS.}

	 In connection with `telescopic products' we shall need later
	that our sections (and homotopies) be locally more regular than
	just plain continuous.
	 Smoothness (say, of class $C^1$) would be quite enough, but is
	hard to achieve if the ground \bspc\ lacks smooth partitions of
	unity.
	 Instead, a notion of boundedness and uniform continuity with respect
	to a preferred trivialization would also do, but it is more convenient
	here
	to use locally Lipschitz sections (written $\llip$), whose definition
	does not require a fixed preferred trivialization.

	 Let $X$ be a separable real \bspc, $\Oa\subset X$ open, $\Ga\to\Oa$ a
	$C^3$-smooth real \blgpbdl, $\dot\Ga\to\Oa$ the associated $C^2$-smooth
	\blabdl\ with a \cts\ norm $|\cdot|_x:\dot\Ga_x\to[0,\infty)$,
	$P\to\Oa$ a $C^3$-smooth right $\Ga$-principal bundle.

	 The main point of this section is the following \t\tHA.

\tetel{\t\tHA.}{For any $f\in C(\Oa,P)$ there are a section
	$g\in\llip(\Oa,P)$ and a homotopy $h\in C(\Oa\times[0,1],P)$
	\st\ $h(x,t)\in P_x$, $h(x,0)=f(x)$, and $h(x,1)=g(x)$ for
	all $(x,t)\in\Oa\times[0,1]$.
}

	 \t\tHA\ follows from \t\tHB\ on approximating a \cts\ section $f$
	by a locally Lipschitz one $g$.

\tetel{\t\tHB.}{For any $f\in C(\Oa,P)$ and $\epsz\in C(\Oa,(0,\infty))$
	there is a $g\in\llip(\Oa,P)$ \st\/
	$|\log_x f(x)^{-1}g(x)|_x\le\epsz(x)$ for all $x\in\Oa$.
}

	 The proof of \t\tHB\ follows that of the Oka principle of Grauert,
	and even if much easier, it still takes a few propositions.

\tetel{\p\tHC.}{Any open covering\/ $\UU$ of an open subset\/ $\Oa$ of $X$ admits
	a Lipschitz partition of unity subordinate to\/ $\UU$.
}

	 Here a Lipschitz partition of unity is one whose members are Lipschitz
	\fns\ $\chi$ on $\Omega$
	the supports of which form a refinement of $\UU$.
	 
\biz{Proof.}{This is a standard fact; see [\rZ].
}

	 We also introduce certain functions $\epsz_i,M_i\in C(\Oa,(0,\infty))$
	which only depend on the geometry of the \blgpbdl\ $\Ga$, such as
	the radius function of a normal \nbd\ of $0$ in $\dot\Ga$.

\tetel{\p\tHD.}{There are $\epsz_1,M_1\in C(\Oa,(0,\infty))$ \st\
	if $U,V\subset\Oa$ open, $h_{UV}\in\llip((U\cap V)\times[0,1],\dot\Ga)$,
	$\la\in\llip((U\cap V)\times[0,1],\End(\dot\Ga))$, $\epsz\in C(\Oa,(0,\infty))$,
	$0<\epsz(x)\le\epsz_1(x)$ for $x\in\Oa$,
	$|h_{UV}(x,t)|_x\le\epsz(x)$,
	$\|\la(x,t)-1\|_x\le\epsz(x)$ for all $(x,t)\in(U\cap V)\times[0,1]$,
	then there are $h_U\in\llip(U\times[0,1],\dot\Ga)$ and
	$h_V\in\llip(V\times[0,1],\dot\Ga)$ \st\
	$h_{UV}(x,t)=h_V(x,t)-\la(x,t)h_U(x,t)$ for $(x,t)\in(U\cap V)\times[0,1]$,
	$|h_U(x,t)|_x\le M_1(x)\epsz(x)$ for $(x,t)\in U\times[0,1]$, and\/
	$|h_V(x,t)|_x\le M_1(x)\epsz(x)$ for $(x,t)\in V\times[0,1]$.
}

	 Here $\|T\|_x$ is the operator norm of a linear operator
	$T\!:\!(\dot\Ga_x,|\cdot|_x)\!\to\!(\dot\Ga_x,|\cdot|_x)$.
	
\biz{Proof.}{Let $0\le\chi_U,\chi_V\le1$, $\chi_U+\chi_V=1$, be a
	locally Lipschitz partition of unity over $U\cup V$ subordinate to
	the covering $\{U,V\}$ and as usual let
	$h_U(x,t)=-\chi_V(x)\la(x,t)^{-1}h_{UV}(x,t)$, and
	$h_V(x,t)=\chi_U(x)h_{UV}(x,t)$.
	 Then $h_{UV}=h_V-\la h_U$, and the rest follows,
	completing the proof of \p\tHD.
}

\tetel{\p\tHE.}{Let $\epsz_0\in C(\Oa,(0,\infty))$ be the radius of a normal
	\nbd\ of\/ $0$ in\/ $\dot\Ga$.
	 Then there are $\epsz_2,M_2\in C(\Oa,(0,\infty))$ with 
	$0<\epsz_2(x)<\epsz_0(x)$ for $x\in\Oa$ \st\ if $x\in\Oa$,
	$0<\epsz\le\epsz_2(x)$, $a,b\in\Ga_x$, $|\log_x a|_x,|\log_x b|_x\le\epsz$,
	then\/ $|\log_x ab|_x\le M_2(x)\epsz$.
}

\tetel{\p\tHF.}{There are $\epsz_3,M_3\in C(\Oa,(0,\infty))$ \st\ if
	$x\in\Oa$, $\xi,\eta\in\dot\Ga_x$, $0<\epsz\le\epsz_3(x)$,
	$|\xi|_x+|\eta|_x\le\epsz$, $\ga\in[0,1]$, then
$$
\eqalign{
	|\log_x[\exp_x(\xi)\exp_x(\ga\log_x(&\exp_x(-\xi)\exp_x(\eta)))]|_x\cr
	&\le
	|(1-\ga)\xi+\ga\eta|_x+M_3(x)(|\xi|_x+|\eta|_x)^2\cr
	&\le\epsz+M_3(x)\epsz^2.\cr
}
	\tag\eHA
 $$
}

\biz{Proofs of Propositions \tHE\ and \tHF.}{These statements follow from the
	second order Taylor formula.
}

	 Let $U\subset\Oa$ be open.
	 Given an $h\in\llip(U\times[0,1],\dot\Ga)$ look at the parametric
	 initial value problem
$$
	\frac{dH}{dt}=H\cdot h,\qquad H(x,0)=1,
\tag\eHB
 $$
 	for an $H\in\llip(U\times[0,1],\Ga)$, where the symbol
	$H(x,t)\cdot h(x,t)$ represents the left translate of the Lie algebra
	element $h(x,t)\in\dot\Ga_x$ by the Lie group element $H(x,t)\in\Ga_x$. 

\tetel{\p\tHG.}{There are functions $\epsz_4,M_4\in C(\Oa,(0,\infty))$ \st\ if
	$\epsz\in C(\Oa,(0,\infty))$ with\/ $0<\epsz(x)\le\epsz_4(x)$ for $x\in\Oa$,
	and $h\in\llip(U\times[0,1],\dot\Ga)$ with\/ $|h(x,t)|_x\le\epsz(x)$
	for $(x,t)\in U\times[0,1]$, then $H\in\llip(U\times[0,1],\Ga)$ as 
	in\/ $(\eHB)$ satisfies that\/
	$|\log_x H(x,t)|_x\le M_4(x)\epsz(x)$ for $(x,t)\in U\times[0,1]$.
}

\biz{Proof.}{This follows from a standard a priori estimate for ordinary differential
	equations.
}

	 Given an $\epsz\in C(\Oa,(0,\infty))$ define $\epsz_5\in C(\Oa,(0,\infty))$ by
$$
	\epsz_5(x)=\frac14\min\{\epsz_0(x),\epsz_1(x),\epsz_2(x),\epsz_4(x),
	\epsz(x)/(1+M_1(x)M_2(x)^2M_4(x))\}.
 $$

\tetel{\p\tHH.}{Let $U,V\subset\Oa$ be open, $f\in C(\Oa,P)$, 
	$\epsz\in C(\Oa,(0,\infty))$.
	 If $g_U\in\llip(U,P)$, $g_V\in\llip(V,P)$ satisfy that\/
	$|\log_x f(x)^{-1}g_U(x)|_x\le\epsz_5(x)$ for $x\in U$,
	and $|\log_x f(x)^{-1}g_V(x)|_x\le\epsz_5(x)$ for $x\in V$, 
	then there is a $g\in\llip(U\cup V,\Ga)$ with\/
	$|\log_x f(x)^{-1}g(x)|_x\le\epsz(x)$ for $x\in U\cup V$.
}

\biz{Proof.}{Define $g_{UV}\in\llip(U\cap V,\Ga)$ by $g_{UV}(x)=g_U(x)^{-1}g_V(x)$.
	 By \p\tHE\ we see that $|\log_x g_{UV}(x)|_x\le M_2(x)\epsz_5(x)$.
	 Define $h_{UV}\in\llip((U\cap V)\times[0,1],\dot\Ga)$ by
	$h_{UV}(x,t)=t\log_x g_{UV}(x)$, $H_{UV}\in\llip((U\cap V)\times[0,1],\Ga)$
	by $H_{UV}(x,t)=\exp_x h_{UV}(x,t)$, and 
	$\la\in\llip((U\cap V)\times[0,1],\End(\dot\Ga))$ by
	$\la(x,t)=\Ad_x(H_{UV}(x,t))$.
	 \p\tHD\ gives $h_U\in\llip(U\times[0,1],\dot\Ga)$ and
	$h_V\in\llip(V\times[0,1],\dot\Ga)$ with
	$|h_U(x,t)|_x\le M_1(x)M_2(x)\epsz_5(x)$ for $(x,t)\in U\times[0,1]$,
	$|h_V(x,t)|_x\le M_1(x)M_2(x)\epsz_5(x)$ for $(x,t)\in V\times[0,1]$, and
$$
	h_{UV}(x,t)=h_V(x,t)-\Ad_x(H_{UV}(x,t))h_U(x,t)
	\tag\eHC
 $$
 	for $(x,t)\in(U\cap V)\times[0,1]$.
	 \p\tHG\ gives $H_U\in\llip(U\times[0,1],\Ga)$ and
	$H_V\in\llip(V\times[0,1],\Ga)$ that satisfy the analogs of $(\eHB)$,
	and the estimates $|\log_x H_U(x,t)|_x\le M_1(x)M_2(x)M_4(x)\epsz_5(x)$ for
	$(x,t)\in U\times[0,1]$, and
	$|\log_x H_V(x,t)|_x\le M_1(x)M_2(x)M_4(x)\epsz_5(x)$ for
	$(x,t)\in V\times[0,1]$.
	 By $(\eHC)$ we see that $H_{UV}(x,t)=H_U(x,t)^{-1}H_V(x,t)$ for
	$(x,t)\in(U\cap V)\times[0,1]$.
	 In particular, for $t=1$ we get that
	$g_U(x)^{-1}g_V(x)=H_U(x,1)^{-1}H_V(x,1)$.
	 Define $g\in\llip((U\cup V)\times[0,1],\Ga)$ by 
	$g(x)=g_U(x)H_U(x,1)^{-1}=g_V(x)H_V(x,1)^{-1}$.
	 Then 
$$\eqalign{
	|\log_x f(x)^{-1}g(x)|_x&=|\log_x f(x)^{-1}g_U(x)H_U(x,1)^{-1}|_x\cr
	&\le M_2(x)M_1(x)M_2(x)M_4(x)\epsz_5(x)<\epsz(x)\cr
}
 $$
	if $x\in U$, and similarly if $x\in V$.
	 The proof of \p\tHH\ is complete.
}

	 Let $f\in C(\Oa,P)$, $U\subset\Oa$ open.  
	 For short we say that $U$
	is \kiemel{good} for $f$ if for any $\epsz\in C(\Oa,(0,\infty))$ there
	is a $g\in\llip(U,P)$ with $|\log_x f(x)^{-1}g(x)|_x\le\epsz(x)$ for
	all $x\in U$.
	 In this language \p\tHH\ says that the union of two (or finitely many)
	good open sets for $f$ is a good open set for $f$.
	 If $U$ is good for $f$, and $V\subset U$ is open, then clearly
	$V$ is also good for $f$.

\tetel{\p\tHI.}{Any point $x_0\in\Oa$ has a good open \nbd\ $U$ for $f$.
}

\biz{Proof.}{Let $U$ be so small an open \nbd\ of $x_0$ that $\Ga$ and $P$ are
	trivial over $U$, and let $T:P|U\to G$ be a trivialization.  
	 Here the \blgp\ $G$ with \bla\ $\dot G$
	is the fiber type of $\Ga$, and define 
	$F:U\to G$ by $F(x)=T(x)f(x)$.
	 It is enough to show that $x_0$ has a good \nbd\ $U$ for $F$.
	 To this end by shrinking $U$ towards $x_0$ we can find a constant
	$c\in G$ \st\ $F(x)c^{-1}$ is so nearly $1$ that $F(x)=c\exp\fii(x)$
	for $x\in U$, where $\fii\in C(U,\dot G)$.
	 Now we only need to show that $U$ is a good \nbd\ for $\fii$, which is
	clear by a Lipschitz partition of unity.
	 The proof of \p\tHI\ is complete.
}

	 Now we claim that there is an exhaustion of $\Oa$ by open sets good for $f$.

\tetel{\p\tHJ.}{There are good open sets\/ $\Oa_n$ for $f$ with
$$
	\Oa_1\subset\overline{\Oa_1}\subset\Oa_2\subset\overline{\Oa_2}\subset\ldots
	\subset\Oa_n\subset\overline{\Oa_n}\subset\ldots,\quad\text{and}\quad
 	\bigcup_{n=1}^\infty\Oa_n=\Oa.
 $$
}

\biz{Proof.}{\p\tHI\ gives a covering $\UU$ of $\Oa$ by good open sets $U$ for $f$,
	and we may suppose that each $U\in\UU$ is bounded and is at a positive
	distance from $X\setminus\Oa$ (in case $\Oa\not=X$).
	 As $\Oa$ is a Lindel\"of space, being a separable metric space, we
	have a countable subcover $V_n$, $n\ge1$, of $\Oa$. 
	 Let $U_n=\bigcup_{i=1}^n V_i$.
	 Define $\Oa_n$, $n\ge1$, by
$
	\Oa_n=\{x\in\Oa:\dist(x,X\setminus U_n)>1/n\}.	
 $
	 Then $\Oa_n$ is an open subset of $\Oa$ and is good for $f$.
	 To check that $\overline\Oa_n\subset\Oa_{n+1}$, let $x\in\overline\Oa_n$.
	 Then $\dist(x,X\setminus U_n)\ge1/n$, and since $U_n\subset U_{n+1}$ we
	 see that $\dist(x,X\setminus U_{n+1})\ge\dist(x,X\setminus U_n)\ge1/n>1/(n+1)$,
	 so $x\in\Oa_{n+1}$.
	  To verify that $\bigcup_{n=1}^\infty\Oa_n=\Oa$, let $x\in\Oa$ be contained
	 in $U_m$ for an $m$.
	  There is an $n\ge m$ for which the ball $B_X(x_0,2/n)$ is contained in 
	 the open set $U_m$.
	  Then $x\in\Oa_n$, since
	 $\dist(x,X\setminus U_n)\ge\dist(x,X\setminus U_m)\ge2/n>1/n$.
	  The proof of \p\tHJ\ is complete.
}

\biz{Proof of \t\tHB.}{Let $\Oa_n$ be as in \p\tHJ, 
	$\epsz_5$ as in \p\tHH, and take a convergent infinite
	product $\prod_{n=1}^\infty a_n=2$ with terms $a_n>1$, 
	e.g., $a_n=(1-(n+1)^{-2})^{-1}$,
	and define an $\epsz_6\in C(\Oa,(0,\infty))$ \st\
	$\epsz_6(x)<\frac14\min\{\epsz_3(x),\epsz_5(x)\}$, and
	$1+4\epsz_6(x)M_3(x)<a_n$ for $x\in\Oa_n\setminus\Oa_{n-1}$, $n\ge2$.

	 We construct a sequence $g_n\in\llip(\Oa_n,P)$ for which
	$g_{n+1}=g_n$ on $\Oa_{n-1}$ for $n\ge2$, and
	$|\log_x f(x)^{-1}g_n(x)|_x\le\epsz_6(x)\prod_{i=1}^{n-1}a_i$
	for $x\in\Oa_n$.
	 If this can be done, then the limit $g$ of $g_n$ as $n\to\infty$
	will do.

	 Since $\Oa_1$ is good for $f$ a $g_1$ can be chosen with
	$|\log_x f(x)^{-1}g_1(x)|_x\le\epsz_6(x)$ for $x\in\Oa_1$,
	and similarly for a $g_2$.

	 Suppose now that $g_1,\ldots,g_n$ for $n\ge2$ have already 
	been defined, and define $g_{n+1}$ as follows.

	 Choose an $h\in\llip(\Oa_{n+1},P)$ with
	$|\log_x f(x)^{-1}h(x)|_x\le\epsz_6(x)$ on $\Oa_{n+1}$,
	which is a good open set for $f$, and a cutoff function
	$\chi\in\llip(\Oa,[0,1])$ with $\chi=1$ on $\overline{\Oa_{n-1}}$,
	and $\chi=0$ on an open \nbd\ of $\overline{\Oa\setminus\Oa_n}$.
	 Define $g_{n+1}$ by $g_{n+1}(x)=h(x)\exp_x(\chi(x)\log_x(h(x)^{-1}g_n(x)))$.

	 Then $g_{n+1}(x)=g_n(x)$ for $x\in\Oa_{n-1}$,
	and $g_{n+1}(x)=h(x)$ for $x\in\Oa_{n+1}\setminus\Oa_n$,
	so the required estimate
	$|\log_x f(x)^{-1}g_{n+1}(x)|_x\le\epsz_6(x)\prod_{i=1}^n a_i$
	holds there.
	 If $x\in\Oa_n\setminus\Oa_{n-1}$, then
$$\eqalign{
	f(x)^{-1}g_{n+1}&(x)=\cr
	f(x&)^{-1}h(x)
	\cdot
	\exp_x(\chi(x)\log_x((f(x)^{-1}h(x))^{-1}\cdot f(x)^{-1}g_n(x))),\cr
}
 $$
 	so by \p\tHF\ we see that
	$|\log_x f(x)^{-1}g_{n+1}(x)|_x\le\epsz_6(x)\prod_{i=1}^n a_i$
	holds there, too.
	 This completes the induction step, and with it the proof of \t\tHB.
}

\biz{Proof of \t\tHA.}{\t\tHB\ gives for $\epsz_0$ as in \p\tHE\ a
	$g\in\llip(\Oa,P)$ with $|\log_x f(x)^{-1}g(x)|_x\le\epsz_0(x)/2$ for
	$x\in\Oa$.
	 Thus letting $h(x,t)=f(x)\exp_x(t\log_x(f(x)^{-1}g(x)))$
	completes the proof of \t\tHA.
}

	 The above proofs easily give approximation by smoother sections
	than locally Lipschitz if the ground \bspc\ admits smoother partitions
	of unity.

	 The general theme of \t\tHA\ is to represent homotopy classes $[f]$
	of \cts\ maps $f:M\to N$ of \bmfd{}s by smoother ones such as locally
	Lipschitz.
	 Similarly the general theme of \t\tHB\ is to approximate \cts\ maps
	$f:M\to N$ of \bmfd{}s by smoother ones such as locally Lipschitz.
	 If $N$ is finite dimensional or embeds in a \bspc\ $Z$ as smooth
	\nbd\ retract, then one can carry out this approximation simpler
	than suggested above for $f\in C(\Oa,G)$ by first approximating
	$f:M\to Z$ by a smoother map $h:M\to Z$ then projecting back onto $N$
	as $g=r\circ h$ with a smooth \nbd\ retraction $r:U\to N$.
	 While this is much shorter when applicable, it is unclear if it
	applies to a map $f:\Oa\to G$ into a general \blgp\ as in \t\tHB.

\alcim{\sI. THE PROOF OF THEOREM \tAC(d).}

	 In this section we complete the proof of \t\tAC(d) that states that a
	\cts{}ly trivial \holo\ \bvbdl\ is \holo{}ally trivial in certain cases.
	 Resume the context and notation of \t\tAC(d).

\biz{Proof of \t\tAC(d).}{Let $G=\GL(Z)$, and $f\in Z^1(\UU,\OO^G)$ be a defining
	cocycle of our \holo\ \bvbdl\ $E\to\Oa$.
	 By \pshdom\ in $\Oa$ there is a Hartogs radius \fn\ $\aa\in\AA$ \st\
	$\BB(\aa)$ is a refinement of $\UU$, and the components $f_{UV}$ of $f$
	are \bu\ on $U\cap V$ for all $U,V\in\BB(\aa)$.
	 As $E$ is \cts{}ly trivial over $\Oa$, there is a \cts\ cochain
	$\fii=(\fii_U)\in C^0(\BB(\aa),C^G)$ with $f_{UV}(x)=\fii_U(x)\fii_V(x)^{-1}$
	for all $x\in U\cap V$, $U,V\in\BB(\aa)$, i.e., $\fii_U(x)=f_{UV}(x)\fii_V(x)$.
	 Thus our $\fii$ can be regarded as a \cts\ global section $\fii$ of a
	\holo\ principal $G$-bundle $P$ over $\Oa$ defined by $f$.
	 \t\tHA\ shows that our \cts\ global section $\fii$ of $P$ is \cts{}ly
	homotopic through global sections of $P$ to a locally Lipschitz global
	section $\psi$ of $P$, i.e., there are $((x,t)\mapsto\fii_U(x,t))\in
	C(U\times[0,1],G)$ with $f_{UV}(x)=\fii_U(x,t)\fii_V(x,t)^{-1}$,
	$\fii_U(x,0)=\fii_U(x)$, and $\psi_U(x)=\fii_U(x,1)$ locally Lipschitz
	for $x\in U$, $U\in\UU$.
	 For a small enough $\ba\in\AA$, $\ba<\aa$, we can write 
	$\psi_U(x)=A_U\exp(\dot\psi_U(x))$ on $U=B_X(x_0,\ba(x_0))$, where $A_U\in G$
	is a constant, and $\dot\psi_U\in\Cbu(U,\dot G)$ is \bu, in fact Lipschitz.
	 Let $\psi_U(x,t)=A_U\exp((1-t)\dot\psi_U(x))$ for $x\in U$, $t\in[0,1]$, $U\in\BB(\ba)$.
	 Define $f'\in\Zs^1(\BB(\ba),\OO^G_1)$ by $f'_{UV}(x,t)=\psi_U(x,t)^{-1}f_{UV}(x)\psi_V(x,t)$.
	 This $f'_{UV}$ is \bu\ on $U\cap V$, $f'_{UV}(x,0)=1$, and $f'_{UV}(x,1)=A_U^{-1} f_{UV}(x)A_V$
	is \holo\ for $x\in U\cap V$, $U,V\in\BB(\ba)$.
	 \t\tFH\ gives a $g'=(g'_U)\in C^0(\BB(\ba),\OO^G_1)$ that resolves $f'$, i.e.,
	$f'_{UV}(x,t)=g'_U(x,t)g'_V(x,t)^{-1}$.
	 Define $g=(g_U)\in C^0(\BB(\ba),\OO^G)$ by $g_U(x)=A_Ug'_U(x,1)$.
	 Then $g$ is \holo\ and satisfies that $f_{UV}(x)=g_U(x)g_V(x)^{-1}$, i.e.,
	our \holo\ \bvbdl\ $E$ is indeed \holo{}ally trivial over $\Oa$.
	 The proof of \t\tAC(d) is complete.
}

\alcim{\sJ. THE PROOF OF THEOREM \tAC(b).}

	 In this section we complete the proof of \t\tAC(b) that states that $E\oplus Z_1$
	is \holo{}ally trivial in certain cases.
	 Resume the context and notation of \t\tAC(b).

\biz{Proof of \t\tAC(b).}{\p\tGA\ tells us that the \holo\ \bvbdl\ $F=E\oplus Z_1\to\Oa$ is \cts{}ly
	trivial.
	 \t\tAC(d) shows then that $F$ is in fact \holo{}ally trivial over $\Oa$.
	 As the fiber type of $F$ is $Z_1\oplus Z_1\cong Z_1$, the proof of \t\tAC(b)
	is complete.
}

\alcim{\sK. THE PROOF OF THEOREM \tAC(c).}

	 In this section we complete the proof of \t\tAC(c) that states that
	$\OO^E$ is acyclic over $\Oa$ in certain cases.
	 Resume the context and notation of \t\tAC(c).

\biz{Proof of \t\tAC(c).}{As $F=E\oplus Z_1\to\Oa$ is \holo{}ally trivial by \t\tAC(b),
	we see by \t\tAC(a) that $0=H^q(\Oa,\OO^F)=H^q(\Oa,\OO^E)\oplus H^q(\Oa,\OO^{Z_1})$
	for $q\ge1$.
	 Thus $H^q(\Oa,\OO^E)=0$ for $q\ge1$, and the proof of \t\tAC(c) is complete.
}

\alcim{\sL. THE PROOF OF THEOREM \tAC(e).}

	 In this section we complete the proof of \t\tAC(e) that says that \holo\
	\hvbdl{}s are \holo{}ally trivial in certain cases.
	 Resume the context and notation of \t\tAC(e), and recall the following
	theorem of Kuiper.

\tetel{\t\tLA.}{{\rm (Kuiper, [\rK])} The \blgp\/ $\GL(\ell_2)$ is contractible,
	and thus any topological \hvbdl\ of fiber type $\ell_2$ over a paracompact
	Hausdorff space (e.g., a metric space) is \cts{}ly trivial.
}

\biz{Proof of \t\tAC(e).}{As $E$ is \cts{}ly trivial over $\Oa$ by \t\tLA\ we can apply
	\t\tAC(d) and find that $E\to\Oa$ is \holo{}ally trivial over $\Oa$ as well.
	The proof of \t\tAC(e) is complete.
}

\alcim{\sX. THE PROOF OF THEOREM \tAC(f).}

	 In this section we finish the proof of \t\tAC(f) that says that a \holo\
	\bvbdl\ $E$ over a contractible space $\Oa$ is \holo{}ally trivial in
	certain cases.
	 Resume the context and notation of \t\tAC(f).

\biz{Proof of \t\tAC(f).}{As $\Oa$ is contractible, our bundle $E\to\Oa$ is \cts{}ly
	trivial by elementary algebraic topology.
	 \t\tAC(d) then shows that $E\to\Oa$ is \holo{}ally trivial.
	 The proof of \t\tAC(f) is complete.
}

\alcim{\sM. APPLICATIONS.}

	 This section points out some applications of \t\tAC.

	 The first application shows that \t\tAC\ remains valid if the open set $\Oa$
	is replaced by certain \cpx\ \bmfd{}s $M$.

\tetel{\t\tMA.}{Let $X$ be a \bspc\ with a \sbs, $\Oa\subset X$ \pscx\ open,
	$M$ a relatively closed \cpx\ Banach submanifold of\/ $\Oa$ onto which there
	is a \holo\ retraction $r:\oa\to M$, where $\oa$ is \pscx\ open with
	$M\subset\oa\subset\Oa$, and $r(x)=x$ if $x\in M$.
	 If \pshdom\ holds in $\oa$, then \t\tAC\ holds with\/ $\Oa$ replaced by $M$.
}

\biz{Proof.}{Consider the pullback bundle $r^*E\to\oa$, and apply to it \t\tAC, and
	then restrict back to $M$.
	 The proof of \t\tMA\ is complete.
}

	\t\tMA\ applies in many cases, e.g., it is shown in [\rPE] that if
	$M$ is a complete intersection in $\Oa$ in the sense that there is a 
	\holo\ \fn\ $f\in\OO(\Oa,Z)$ into a \bspc\ $Z$, $M=\{x\in\Oa:f(x)=0\}$,
	and the Fr\'echet differential $df(x)\in\Hom(X,Z)$ has split kernel for
	$x\in M$, then there is a \holo\ retraction $r:\oa\to M$ as in \t\tMA.

	 The total space of a direct summand of a trivial \holo\ \bvbdl\ is a
	complete intersection in the trivial bundle, hence is the following theorem
	that proves a vanishing result over the total space of a \holo\ \bvbdl\ in
	certain cases.

\tetel{\t\tMX.}{Let $X$ be a \bspc\ with a \sbs, $\Oa\subset X$ \pscx\ open,
	$E'\to\Oa$ a \holo\ \bvbdl\ with a \bspc\ $Z'$ for fiber type,
	$E\to E'$ a \holo\ \bvbdl\ over the total space of the bundle $E'\to\Oa$,
	and $1\le p<\infty$.
	 If $Z'_1=\ell_p(Z')$ has a \sbs, and \pshdom\ holds in\/ $\Oa\times Z'_1$,
	then \t\tAC\ holds with\/ $\Oa$ replaced by $E'$.
}

\biz{Proof.}{\t\tAC(b) allows us to regard $E'$ as a direct summand in the trivial
	bundle $\Oa\times Z'_1$, and provides a projection $P\in\OO(\Oa,\End(Z'_1))$
	with $P(x)^2=P(x)$, and $\im\, P(x)=E'_x$, $x\in\Oa$.
	 Thus $r:\Oa\times Z'_1\to E'\subset\Oa\times Z'_1$ defined by
	$r(x,\xi)=(x,P(x)\xi)$ is a \holo\ retraction, and an application of \t\tMA\
	completes the proof of \t\tMX.
}

\tetel{\t\tMB.}{If $X=\ell_2$ in \t\tMA\, and $M$ is \idml, then the \holo\ tangent bundle $TM$
	and cotangent bundle $T^*M$ of $M$ are \holo{}ally trivial over $M$, also 
	so are the bundles $S^pM$, $\La^pM$ of \holo\ symmetric or alternating
	$p$-forms on $M$ for $p\ge1$.
}

\biz{Proof.}{As $TM$ and $T^*M$ are Hilbert bundles of fiber type $\ell_2$, 
	they are \holo{}ally trivial by \t\tMA. As the bundles $S^pM,\La^pM$ are
	associated bundles to the trivial bundle $TM$, they are themselves \holo{}ally
	trivial---they are defined by an action of the defining cocycle $f$ of $TM$ on
	the \bspc{}s $S^p\ell_2,\La^p\ell_2$, but as the defining cocycle $f$ is
	\holo{}ally trivial, so are the associated bundles $S^pM,\La^pM$.
	The proof of \t\tMB\ is complete.
}

	 So for example if $M$ is as in \t\tMB,
	 then there are a \holo\ $(1,0)$-form $\aa$ and a \holo\ $(1,0)$-vector field $\xi$ on $M$ 
	that have no zeros on $M$.
	 Can such an $\aa$ be chosen exact on $M$, i.e., of the form $\aa=df$, where
	$f\in\OO(M)$?
	 Can such a $\xi$ be chosen complete for real time on $M$?

	 We now look at analytic subsets $A$ of $\Oa$ and their \nbd\ bases.

\tetel{\t\tMC.}{With the notation and hypotheses of \t\tAC, any analytic subset $A$
	of\/ $\Oa$ that can be defined as a set by $A=\{x\in\Oa:s(x)=0\}$, where
	$s\in\OO(\Oa,E)$ is a \holo\ section of a \holo\ \bvbdl\ $E\to\Oa$ can in fact be
	defined as a set by $A=\{x\in\Oa:f(x)=0\}$, where $f\in\OO(\Oa,Z_1)$ is a
	\holo\ \fn\ with values in a \bspc\ $Z_1$.
}

\biz{Proof.}{\t\tAC(b) gives an injective map $I\in\OO(\Oa,\Hom(E,Z_1))$ of \holo\
	\bvbdl{}s.
	 Defining $f\in\OO(\Oa,Z_1)$ by $f(x)=I(x)s(x)$ completes the proof of \t\tMC.
}

	 \t\tMC\ applies in many cases.
	 If $A\subset\Oa$ is a possibly singular hypersurface,
	or an iterated hypersurface in the sense that $A=M_n\subset M_{n-1}\subset M_1\subset M_0=\Oa$,
	where $M_i$ is smooth \cpx\ hypersurface in $M_{i-1}$, $i=1,\ldots,n$, then
	$A$ can be defined as a set by $A=\{x\in\Oa':s(x)=0\}$ as in \t\tMC,
	where $\Oa'$ is \pscx\ open with $A\subset\Oa'\subset\Oa$.
	 See [\rPE], where it is also shown that if $A$ can be defined as a set
	by $A=\{x\in\Oa:s(x)=0\}$ as in \t\tMC, then $A$ has a \nbd\ basis consisting
	of \pscx\ open subsets of $\Oa$.

	 In classical Stein theory the sheaves $\OO^Z$, $Z=\CC^n$, $n\ge1$, are called
	free sheaves and serve as building blocks in acyclic resolutions of more general
	analytic sheaves, called coherent analytic sheaves.
	 Let us call the sheaf $\OO^Z$ a model sheaf, where $Z$ is any \bspc.
	 These model sheaves can be used to define a class of sheaves that is
	analogous to the class of coherent analytic sheaves in \fdms, see [\rPC].

\tetel{\t\tMD.}{With the notation and hypotheses of \t\tAC.
\vskip0pt
	{\rm(a)} There is a \holo\ operator-valued \fn\ $A\in\OO(\Oa,\End(Z_1))$ \st\
	$E$ is \holo{}ally isomorphic over\/ $\Oa$ to the cokernel of $A$, i.e.,
	there is a locally (and in fact globally) split \ses\/
	$0\to\Oa\times Z_1\to\Oa\times Z_1\to E\to0$
	of \holo\ \bvbdl{}s over\/ $\Oa$ with the second map being
	$(x,\xi)\mapsto(x,A(x)\xi)$
\vskip0pt
	{\rm(b)} The sheaf $\OO^E$ has a global resolution\/
	$0\to\OO^{Z_1}\to\OO^{Z_1}\to\OO^E\to0$ over\/ $\Oa$ by model sheaves
	that is exact on the level of germs and on the level of global section over
	any \pscx\ open subset $U$ of\/ $\Oa$, if \pshdom\ holds in $U$.
}

\biz{Proof.}{Part (a) follows easily from \t\tAC(b) while (b) follows from (a)
	noting that any locally split \ses\ $0\to E'\to E\to E''\to0$ 
	of \holo\ \bvbdl{}s is in fact globally split since $H^1(\Oa,\OO^F)=0$
	for $F=\Hom(E'',E')$ by \t\tAC(c).
	 The proof of \t\tMD\ is complete.
}

	 Vector bundles are also applicable to some questions on the algebra of
	\holo\ \fns\ such as the weak Nullstellensatz or completion of \holo\
	matrices to invertible \holo\ matrices.

\tetel{\t\tME.}{Let $X$ be a \bspc\ with a \sbs, $\Oa\subset X$ \pscx\ open,
	and $f=(f_n)\in\OO(\Oa,\ell_2)$ with $f(x)\not=0$ for $x\in\Oa$.
	 Then there is a $g=(g_n)\in\OO(\Oa,\ell_2)$ with\/
	$\sum_{n=1}^\infty f_n(x)g_n(x)=1$ for $x\in\Oa$.
}

\biz{Proof.}{Look at the \ses\ $0\to K\to\Oa\times\ell_2\to\Oa\times\CC\to0$
	of \holo\ \hvbdl{}s over $\Oa$, where the second map is inclusion and
	the third is $(x,\xi)\mapsto F(x,\xi)=\sum_{n=1}^\infty f_n(x)\xi_n$.
	 As $H^1(\Oa,\OO^K)=0$ by \t\tAC(c), our \ses\ splits, i.e., there is
	a map $\xi=g(x)$ with $F(x,g(x))=1$ for $x\in\Oa$.
	 The proof of \p\tME\ is complete.
}

\tetel{\t\tMF.}{Let $X$ be a \bspc\ with a \sbs, $\Oa\subset X$ \pscx\ open,
	and $f_1,\ldots,f_n\in\OO(\Oa,\ell_2)$, $n\ge1$, column vectors that
	are pointwise linearly independent in $\ell_2$.
	 If \pshdom\ holds in\/ $\Oa$, then there are further columns
	$f_{n+1},f_{n+2},\ldots\in\OO(\Oa,\ell_2)$ with the property that
	the matrix $A(x)$ whose $(i,j)$ component is the $i$'th component 
	of the vector $f_j(x)$ defines an invertible operator $A\in\OO(\Oa,\GL(\ell_2))$
	the first $n$ columns of which are the given $f_1,\ldots,f_n$.
}

\biz{Proof.}{Look at the \ses\ $0\to\Oa\times\CC^n\to\Oa\times\ell_2\to E\to0$
	of \holo\ \hvbdl{}s, where the second map is $(x,\xi)\mapsto\sum_{i=1}^n f_i(x)\xi_i$.
	 As $E$ is \holo{}ally isomorphic to $\Oa\times\ell_2$ by \t\tAC(d), and
	as the above \ses\ \holo{}ally splits over $\Oa$ by \t\tAC(c), the proof of \t\tMF\ is complete.
}

	 Below is a version of \t\tMF\ for a finite matrix of determinant $1$.

\tetel{\t\tMG.}{Let $X$ be a \bspc\ with a \sbs, $\Oa\subset X$ \pscx\ open,
	and $f_1,\ldots,f_n\in\OO(\Oa,\CC^N)$, $1\le n<N$, column vectors that
	are pointwise linearly independent in\/ $\CC^N$.
	 If \pshdom\ holds in\/ $\Oa$, and there are further \cts\ columns
	$f_{n+1},\ldots,f_N\in C(\Oa,\CC^N)$ with the property that
	the matrix $A(x)$ whose $(i,j)$ component is the $i$'th component
	of the vector $f_j(x)$ 
	has determinant\/ $1$ for $x\in\Oa$, and the first $n$ columns 
	of which are the given $f_1,\ldots,f_n$,
	then there are \holo\ $f_{n+1},\ldots,f_N\in\OO(\Oa,\CC^N)$ with the
	same property.
}

\biz{Proof.}{Look at the \ses\ $0\to\Oa\times\CC^n\to\Oa\times\CC^N\to E\to0$
	of \holo\ \vbdl{}s, where the second map is $(x,\xi)\mapsto\sum_{i=1}^n f_i(x)\xi_i$.
	 As $E$ is \cts{}ly trivial by the assumption about the existence
	of \cts\ augmentation $f_{n+1},\ldots,f_N$, we see by \t\tAC(d) that
	$E$ is \holo{}ally trivial as well over $\Oa$, and as the above \ses\
	\holo{}ally splits over $\Oa$ by \t\tAC(c), the proof of \t\tMG\ is complete.
	Note that if $\Oa$ is contractible, say, star-like or convex, then
	there is a continuous augmentation, and \t\tMG\ applies.
}

	 The meaning of \t\tMG\ is that under some conditions one (or several
	independent) unimodular row(s) of the ring $\OO(\Oa)$ 
	is (are) equivalent to the standard
	unimodular row $(1,0,\ldots,0)$ (or rows 
	$(1,0,\ldots,0)$, $(0,1,0,\ldots,0)$, \dots, $(0,\ldots,0,1,0,\ldots,0)$)
	in the sense of algebra.

	 We now look at approximation and interpolation for \holo\ sections
	of a \holo{} \bvbdl\ in certain cases.

\tetel{\t\tMH.}{With the notation and hypotheses of \t\tAC\ the following version
	of Runge approximation holds for \holo\ sections of $E$.
	 There are an embedding $I\in\OO(\Oa,\Hom(E,Z_1))$ of $E$ into the trivial
	bundle\/ $\Oa\times Z_1$ as a direct summand $IE$, and a projection
	$P\in\OO(\Oa,\End(Z_1))$ onto $IE$.
	 There is an $\aa$ with $8\aa\in\AA$ with $P$ \bu\ on\/ $\Oa_N\<\aa\>$
	for all $N\ge1$.
	 Choose $\ga\in\AA$ as in \p\tBB(b), and $\ga'\in\AA$ as in \p\tBB(a).
	 Let $\epsz>0$, $N\ge1$, and $f\in\OO(\Oa_N\<\aa\>,E)$ \st\ $If\in\OO(\Oa_N\<\aa\>,Z_1)$
	is \bu.
	 Then the following hold.
\vskip0pt
	{\rm(a)} There is a $g\in\OO(\Oa_{N+1}\<\aa\>,E)$ with $Ig$ \bu\ \st\/
	$\|I(x)(f(x)-g(x))\|<\epsz$ for $x\in\Oa_N\<\ga\>$.
\vskip0pt
	{\rm(b)} There is a $g\in\OO(\Oa,E)$ with $Ig$ \bu\ on\/ $\Oa_{N+p}\<\ga'\>$ for all
	$p\ge0$ \st\/ $\|I(x)(f(x)-g(x)\|<\epsz$ for $x\in\Oa_N\<\ga'\>$.
}

\biz{Proof.}{\t\tAC(b) provides $I$ and $P$, and for a given $P$ \pshdom\ in $\Oa$
	a required Hartogs radius \fn\ $\aa$.
	 An application of \p\tCC\ completes the proof of \t\tMH.
}

\tetel{\t\tMI.}{With the notation and hypotheses of \t\tAC.
	 Let $x_n$, $n\ge1$, be a sequence of points of\/ $\Oa$, $\epsz_n>0$, and
	suppose that the balls $B_n=B_X(x_n,2\epsz_n)$ are disjoint from one another,
	and that the set\/ $\{x_n:n\ge1\}$ has no limit points in\/ $\Oa$.
	 Let $m_n\ge0$ be integers for $n\ge1$.
\vskip0pt
	{\rm(a)} Let $f_n\in\OO(B_n,Z)$ for $n\ge1$.
	 Then there is an $f\in\OO(\Oa,Z)$ with $\|f(x)-f_n(x)\|=\bigO(\|x-x_n\|^{m_n+1})$
	as $x\to x_n$ for all $n\ge1$.
	 Let $J\to\Oa$ be the sheaf of germs of \holo\ \fns\/ $\Oa\to Z$ that vanish 
	at the points $x_n$ to order $m_n$ for all $n\ge1$.
	 Then $H^q(\Oa,J)=0$ for all $q\ge1$.
\vskip0pt
	{\rm(b)} Suppose that $\epsz_n$, $n\ge1$, are so small that $E|B_n$ has a trivialization
	$T_n\in\OO(B_n,\Hom(E,Z))$, and let $f_n\in\OO(B_n,E)$ be local sections.
	 Then there is an $f\in\OO(\Oa,E)$ with\/
	$\|T_n(x)(f(x)-f_n(x))\|=\bigO(\|x-x_n\|^{m_n+1})$ as $x\to x_n$ for all $n\ge1$.
	 Let $J$ be the sheaf of germs of \holo\ sections\/ $\Oa\to E$ that vanish
	at the points $x_n$ to order $m_n$ for all $n\ge1$.
	 Then $H^q(\Oa,J)=0$ for all $q\ge1$.
}

\biz{Proof.}{In both (a) and (b) the vanishing $H^q(\Oa,J)=0$ for $q\ge1$ follows from
	\t\tAC(ac) in a standard way if the said interpolations are possible.

	 To prove interpolation in (a) choose by \pshdom\ in $\Oa$ an $\aa$ with
	$10\aa\in\AA$ so that the covering $\BB(\aa)$ by radius $\aa$ balls is a
	refinement of the covering $\{\Oa\setminus\bigcup_{n=1}^\infty\overline{B_n}, B_n:n\ge1\}$
	of $\Oa$.
	 Then in the exhausting sets $\Oa_N\<\aa\>$ only finitely many points $x_n$ may
	fall for all $N\ge1$.
	 Complete the proof as in [\rPF, \S\,8] relying on \p\tCB\ here for Runge
	approximation there.

	 To prove interpolation in (b) choose $I$ and $P$ as in the proof of \t\tMH\
	and apply (a) to interpolate $If_n$ by an $f'$, and project back to define
	$f$ by $f=Pf'$.
	 The proof of \t\tMI\ is complete.
}

	 Note that in Theorems~\tMC--\tMI\ the open set $\Oa$ can be replaced by a
	\cpx\ \bmfd\ $M$ as in \t\tMA.

\alcim{\sN. DISCUSSION.}

	 In this section we make some remarks about the methods that we chose
	in this paper.

	 In the vanishing theorems of classical Stein theory such as Theorem~B
	and the Grauert--Oka principle the main ingredient are \holo{}ally convex
	compact exhaustions, Runge-type approximation for \holo{} \fns\ over
	\holo{}ally convex compact sets, the Oka coherence theorem, and
	some generalities on sheaf theory.
	
	 The Oka coherence theorem does not seem to have an analog in \idms,
	but since most of the interesting bundles such as the tangent bundle
	and its associated bundles are of infinite rank, the Oka coherence
	theorem, which is in essence a finiteness result, could not be profitable.
	 It is an interesting question to find a class of analytic sheaves
	for which an analog of the classical Theorem~B can be proved.

	 The most important properties of compact sets for the purposes of
	vanishing theorems seem to be that over them \cts\ and \holo\ \fn{}s
	are \bu.  While boundedness would have sufficed in the proof of \t\tAC(a)
	for the additive case of $Z$, the uniform continuity was crucial
	for the classical trick of telescopic products in connection with
	$\GL(Z)$ to go through.  It is true that over sets $\Oa(D,R)$ as in
	$(\eCA)$ uniform continuity over an arbitrarily slightly smaller set
	of the same type can be deduced from boundedness by the Schwarz lemma for \holo\
	\fns, it was much more convenient to build it in the definitions.
	 Especially so for \t\tFH\ since \cts\ \fns\ in \idms\ are
	not necessarily locally \ucts.
	 Another important property of compact sets is that they can be
	covered by finitely many open sets at an arbitrary scale.
	 This, while cannot be fully matched in infinite dimensions with
	open sets of exhaustion, can be modelled by considering a fixed
	scale $\aa\in\AA$ as a Hartogs radius \fn\ at which the object
	at hand (say, a cocycle to be resolved) can be handled.

	 The Runge approximation was necessary only for \bu\ \fns\
	with values in a \bspc\ $Z$ or a \blgp\ $G=\GL(Z)$
	for which it is straightforward.

	 The exhaustion was made possible by \pshdom\ which in essence
	seems to be a useful \holo\ convexity property that allows one
	to make flexible exhaustions to match the size properties of the
	data.  In the classical case one Stein exhaustion works for
	all cocycles of a coherent analytic sheaf, but in \idms\ 
	such a one size fits all approach seems unavailable.
	 The greatest virtue of \pshdom\ is that it allows the construction
	of a fine scale based on the data at which to work.
	 
	 It was L\'aszl\'o Lempert who introduced in [\rLA] 
	(after some precursors; see [\rN]) the way to exhaust a \pscx\
	open set by so-called open sets of type~(B) as in $(\eCA)$ here
	along with the fairly classical way of handling \fns\ over one
	open set of type~(B), in [\rLB] the idea of \pshdom, and in [\rLC]
	how to use it to make an exhaustion tailored to the data.

	 The methods of this paper build on all of these ideas of Lempert
	along with some useful changes and twists.
	 Our methods have the advantage (besides proving stronger results
	than in [\rLC]) that they make up a last on which a number
	of other vanishing theorems can be fashioned.
	 Notice the analogous proofs of \t\tAC(a) and \t\tFH.
	 Further applications include the proof of a fuller version of the
	Grauert--Oka principle in [\rPD], and amalgamation of syzygies in [\rPC].

	 \sbs\ for the ground \bspc\ $X$ is not fully necessary in \t\tAC, since
	the requirement that $X$ have a \sbs\ can be replaced by demanding
	that $X$ be a direct summand in a \bspc\ with a \sbs, i.e., $X$ have
	the bounded approximation property by a theorem of Pe\l czy\'nski's.
	 Regarding the other hypothesis of \t\tAC\ that in $\Oa$ \pshdom\ hold,
	it is reasonable to hope (while not currently proved) that it may
	hold in any \pscx\ open $\Oa$ subset of a \bspc\ $X$ with the bounded
	approximation property.
	 Another advantage of \pshdom\ over Runge approximation on balls seems
	to be that it can be formulated in \bmfd{}s, it seems more directly
	relevant to the business of vanishing, and also it appears to inherit
	better to submanifolds.

	{\vekony Acknowledgements.} The work on this paper started at the Riverside
	campus of the University of California, continued at its San Diego
	campus, and was completed at Georgia State University.
	 The author is grateful to these institutions and to 
	Professors Bun Wong of UCR, Peter Ebenfelt of
	UCSD, and Mih\'aly Bakonyi of GSU, whose support has been essential to him.
	The author also thanks Professor L\'aszl\'o Lempert for his comments
	on the manuscript of this paper.

\vskip0.30truein
\centerline{\scVIII References}
\vskip0.20truein
\baselineskip=11pt
\parskip=7pt
\frenchspacing
{\rmVIII

	[\rC] Cartan,~H., 
	{\itVIII 
	Espaces fibr\'es analytiques},
	Symposium Internacional de Topologia Algebraica, Mexico, 
	(1958), 97--121.
	
	[\rZ] Frol\'\i k,~Z.,
	{\itVIII
	Existence of $\scriptstyle\ell_\infty$-partitions of unity},
	Rend. Sem. Mat. Univ. Politec. Torino,
	{\bfVIII 42}
	(1985),
	9--14.

	[\rH] H\"ormander,~L., 
	{\itVIII
	An Introduction to Complex Analysis in Several Variables},
	3rd Ed., North-Holland, Amsterdam,
	(1990).

	[\rK] Kuiper,~N.H.,
	{\itVIII
	The homotopy type of the unitary group of Hilbert space},
	Topology,
	{\bfVIII 3}
	(1965),
	19--30.
	
	[\rLt] Leiterer,~J.,
	{\itVIII
	Banach coherent analytic Fr\'echet sheaves},
	Math.{} Nachr., {\bfVIII 85} (1978), 91--109.

	[\rLA] Lempert,~L.,
    	{\itVIII
   	The Dolbeault complex in infinite dimensions~III},
    	Invent.{} Math., {\bfVIII 142} (2000), 579--603.

	[\rLB] \vonal,
	{\itVIII
	Plurisubharmonic domination},
	J. Amer. Math. Soc.,
	{\bfVIII 17}
	(2004),
	361--372.

	[\rLC] \vonal,
	{\itVIII
	Vanishing cohomology for holomorphic vector bundles in a Banach setting},
	Asian J. Math.,
	{\bfVIII 8}
	(2004),
	65--85.

    	[\rLD] \vonal,
    	{\itVIII
	Acyclic sheaves in Banach spaces},
	Contemporary Math.,
	{\bfVIII 368}
	(2005),
	313--320.

	[\rN] Noverraz,~P.,
	{\itVIII
	Pseudo-convexit\'e, convexit\'e polynomiale et domains
	d'holomorphie en dimension infinie},
	North--Holland, Amsterdam, (1973).
	
	[\rPA] Patyi,~I.,
	{\itVIII
	On the Oka principle in a Banach space I},
	Math. Ann.,
	{\bfVIII 326}
	(2003),
	417--441.

	[\rPF] \vonal,
	{\itVIII
	Cohomological characterization of pseudoconvexity in a Banach space},
	Math. Z.,
	{\bfVIII 245}
	(2003),
	371--386.

	[\rPE] \vonal,
	{\itVIII
	Analytic cohomology of complete intersections in a Banach space},
	Ann. Inst. Fourier (Grenoble),
	{\bfVIII 54}
	(2004),
	147--158.

	[\rPB] \vonal,
	{\itVIII
	An analytic Koszul complex in a Banach space},
	manuscript, (2005).

	[\rPC] \vonal,
	{\itVIII
	Analytic cohomology in a Banach space},
	in preparation.

	[\rPD] \vonal,
	{\itVIII
	The Grauert--Oka principle in a Banach space},
	in preparation.

\vskip0.20truein
\centerline{\vastag*~***~*}
\vskip0.15truein
{\scVIII
	Imre Patyi,
	Department of Mathematics and Statistics,
	Georgia State University,
	Atlanta, GA 30303-3083, USA,
	{\ttVIII ipatyi\@gsu.edu}}
\bye